\documentclass[final,twocolumn,5p,times,10pt]{elsarticle}
\usepackage{lineno}
\modulolinenumbers[5]

\journal{Automatica}

\usepackage[utf8]{inputenc}
\usepackage[american]{babel}

\usepackage{amsthm}
\usepackage{amsmath}
\usepackage{amssymb}
\usepackage{amsfonts}
\usepackage{bm}
\usepackage{hhline}
\usepackage{csquotes}
\usepackage{color}
\usepackage{import}
\usepackage{algorithm}
\usepackage{algpseudocode}
\usepackage{dsfont}
\usepackage{multirow}
\usepackage{tensor}
\usepackage{subfig}
\usepackage{soul}
\usepackage{graphicx}
\usepackage{microtype}
\usepackage{dsfont}

\usepackage{pdfpages}

\makeatletter
\usepackage[colorlinks=true,linkcolor=red]{hyperref}
\makeatother

\newtheorem{theorem}{Theorem}
\newtheorem{lemma}{Lemma}
\newtheorem{proposition}{Proposition}
\newtheorem{remark}{Remark}
\newtheorem{definition}{Definition}
\newtheorem{property}{Property}

\allowdisplaybreaks 

\newcommand{\gzinclusion}{\ensuremath{\triangleleft}}

\newcommand{\ninf}[1]{\ensuremath{\|{#1}\|_\infty}}
\newcommand{\none}[1]{\ensuremath{\|{#1}\|_1}}
\newcommand{\conball}{\ensuremath{B_\infty(\mathbf{A},\mathbf{b})}}

\newcommand{\xim}{\ensuremath{\bm{\xi}_\text{m}}}

\newcommand{\xil}{\ensuremath{\bm{\xi}^\text{L}}}
\newcommand{\xiu}{\ensuremath{\bm{\xi}^\text{U}}}

\newcommand{\arrowmatrix}[3]{\ensuremath{\underset{\xrightarrow[\scriptstyle{#2}]{\hphantom{#1}}}{\left[#1\right]}\left.\!\!\vphantom{#1}\right\downarrow\!{\scriptstyle{#3}}}}

\newcommand{\downarrowmatrix}[2]{\ensuremath{{\left[#1\right]}\left.\!\!\vphantom{#1}\right\downarrow\!{\scriptstyle{#2}}}}

\newcommand{\textred}[1]{\textcolor{red}{#1}}

\newtheorem{corollary}{\bf{Corollary}}

\newcommand{\half}{\ensuremath{\frac{1}{2}}}

\newcommand{\eyenoarg}{\ensuremath{\mathbf{I}}}
\newcommand{\zeros}[2]{\ensuremath{\bm{0}_{#1\times#2}}}

\newcommand{\realset}{\ensuremath{\mathbb{R}}}
\newcommand{\realsetmat}[2]{\ensuremath{\mathbb{R}^{#1\times#2}}}

\newcommand{\trace}[1]{\ensuremath{\text{trace}\left\{{#1}\right\}}}

\newcommand{\lbound}{\ensuremath{\text{L}}}
\newcommand{\ubound}{\ensuremath{\text{U}}}
\newcommand{\midpoint}[1]{\ensuremath{\text{mid}({#1})}}

\newcommand{\diam}[1]{\ensuremath{\text{diam}({#1})}}
\newcommand{\rad}[1]{\ensuremath{\text{rad}({#1})}}

\newcommand{\intvalset}{\ensuremath{\mathbb{IR}}}
\newcommand{\intvalsetmat}[2]{\ensuremath{\intvalset^{{#1} \times {#2}}}}

\newcommand{\iextension}[1]{\ensuremath{\square \left( #1 \right) }}

\newcommand{\alamobravo}{ZMV}
\newcommand{\combastelbravo}{ZFO}

\newcommand{\mbf}[1]{\ensuremath{\mathbf{#1}}}

\begin{document}

\begin{frontmatter}

\title{Guaranteed methods based on constrained zonotopes for set-valued state estimation of nonlinear discrete-time systems\tnoteref{mytitlenote}}
\tnotetext[mytitlenote]{{\copyright} 2019. This manuscript version is made available under the CC-BY-NC-ND 4.0 license \url{https://creativecommons.org/licenses/by-nc-nd/4.0/}. This work was supported in part by the project INCT under the grant CNPq 465755/2014-3, FAPESP 2014/50851-0, and also by the Brazilian agencies CAPES, and FAPEMIG. The material in this paper was not presented at any conference.}



\author[myfirstaddress]{Brenner S. Rego\corref{mycorrespondingauthor}}
\cortext[mycorrespondingauthor]{Corresponding author}
\ead{brennersr7@ufmg.br}

\author[myfirstaddress,mysecondaddress]{Guilherme V. Raffo}
\ead{raffo@ufmg.br}

\author[mythirdaddress]{Joseph K. Scott}
\ead{jks9@clemson.edu}

\author[myfourthaddress]{Davide M. Raimondo}
\ead{davide.raimondo@unipv.it}

\address[myfirstaddress]{Graduate Program in Electrical Engineering, Federal University of Minas Gerais, Belo Horizonte, MG 31270-901, Brazil}
\address[mysecondaddress]{Department of Electronics Engineering, Federal University of Minas Gerais, Belo Horizonte, MG 31270-901, Brazil}
\address[mythirdaddress]{Department of Chemical and Biomolecular Engineering, Clemson University, Clemson, SC, USA}
\address[myfourthaddress]{Department of Electrical, Computer and Biomedical Engineering, University of Pavia, Italy}

\begin{abstract}
This paper presents new methods for set-valued state estimation of nonlinear discrete-time systems with unknown-but-bounded uncertainties. A single time step involves propagating an enclosure of the system states through the nonlinear dynamics (prediction), and then enclosing the intersection of this set with a bounded-error measurement (update). When these enclosures are represented by simple sets such as intervals, ellipsoids, parallelotopes, and zonotopes, certain set operations can be very conservative. Yet, using general convex polytopes is much more computationally demanding. To address this, this paper presents two new methods, a mean value extension and a first-order Taylor extension, for efficiently propagating constrained zonotopes through nonlinear mappings. These extend existing methods for zonotopes in a consistent way. Examples show that these extensions yield tighter prediction enclosures than zonotopic estimation methods, while largely retaining the computational benefits of zonotopes. Moreover, they enable tighter update enclosures because constrained zonotopes can represent intersections much more accurately than zonotopes.
\end{abstract}


\begin{keyword}
Nonlinear state estimation, Set-based computing, Reachability analysis, Convex polytopes
\end{keyword}

\end{frontmatter}


\section{Introduction} \label{sec:introduction}


The importance of state estimation has become evident in many fields of research over the years. Examples of recurrent applications are the localization problem \cite{Jaulin2009b,Saeedi2016}, state-feedback control \cite{Jaulin2009,Santos2017,Goodarzi2017}, and fault detection and isolation (FDI) \cite{Combastel2008,Zhao2015,Raimondo2016}. If stochastic descriptions of the process and measurement uncertainties are available, then Bayesian state estimation methods such as Kalman filtering or particle filtering are typically employed. On the other hand, if only bounds on the uncertainties are known, then set-valued state estimation methods are applied \cite{Jaulin2009,Chisci1996,Le2013,Rego2018b,Rego2019a}. Recent approaches also consider the presence of both kinds of uncertainties in the system \cite{Combastel2015}.

Set-valued state estimation methods aim to construct compact sets that are guaranteed to enclose all possible trajectories of the system subject to unknown-but-bounded uncertainties \cite{Chisci1996,Alamo2005a}. Using the standard recursive approach for discrete-time systems, this involves first bounding the image of the current enclosure under the dynamics (prediction), and then enclosing the intersection of this set with the set of states consistent with a bounded-error measurement (update). For discrete-time linear systems, if the initial set is a polytope, then exact enclosures can in principle be computed using convex polytopes \citep{Girard2008}. However, even for linear dynamics, polytope propagation requires demanding operations (e.g., polytope projection, Minkowski sum, or conversion between vertex and halfspace representations) whose complexity grows dramatically with time due to the complexity increase of the resulting sets \cite{Walter1989,Shamma1999}. For these reasons, enclosures are often described by simpler sets including ellipsoids \cite{Schweppe1968,Durieu2001,Polyak2004}, parallelotopes \cite{Chisci1996,Vicino1996}, zonotopes \cite{Le2013,Combastel2003}, or combinations of these \cite{Chabane2014}. However, the mathematical limitations of these sets require certain operations to be conservative, sometimes quite significantly. Notably, this includes set intersection, which is critical for the update step in set-valued state estimation \cite{Chisci1996,Le2013,Durieu2001}. The article \cite{Althoff2011} proposes the use of zonotope bundles to describe intersections of zonotopes without explicit computation. However, the Minkowski sum and linear image (see Section \ref{sec:preliminaries}) are outer-approximated. In the article \cite{Scott2016}, \emph{constrained zonotopes} are introduced to overcome many of the limitations of zonotopes. These sets are closed under intersection, Minkowski sum, and linear image, and are capable of describing arbitrary convex polytopes if the complexity of the set description is not limited. Efficient algorithms for linear set-valued state estimation and also FDI using constrained zonotopes are described in \cite{Scott2016,Raimondo2016}.

In contrast to the linear case, effective set-valued state estimation for nonlinear systems is still an open challenge \cite{Jaulin2009,Alamo2005a,Combastel2005,Wan2018}. Early approaches in this field used inclusion functions based on interval arithmetic \cite{Moore2009} to propagate bounds through the nonlinear dynamics, and used interval-based set inversion techniques to enclose the set of states consistent with the current measurement \cite{Jaulin2009b,Jaulin2009,Kieffer1998}. Improved accuracy is achieved using refinements (i.e., unions of intervals) as well as more advanced interval methods such as contractor and separator algebras \cite{Jaulin2009b,Jaulin2016}. Unfortunately, even these methods often provide conservative bounds without extensive refinement, which is only tractable for systems with relatively few states \cite{Jaulin2009}. An interesting new interval method based on discrete-time differential inequalities is proposed in \cite{Yang2018}, but it only applies to Euler discretized systems with limited step size.

A few alternatives for nonlinear set-valued state estimation can be found in the literature. Polytopes are used in \cite{Shamma1997} for systems with nonlinear dynamics and linear measurements. The prediction step is performed based on a linearization of the dynamics, in which conservative interval enclosures are used to bound the linearization error. Notably, the computed error bound is valid for any value of the state, rather than being computed as a function of the current state enclosure, which can lead to a very conservative result. Another issue with this method is that the complexity of the polytopic enclosures grows rapidly with time, which results in a very high computational burden because the complexity of the required set operations scales poorly with increasing polytope complexity. More efficient methods based on zonotopes are proposed in \cite{Alamo2005a} and \cite{Combastel2005}. The propagation step in \cite{Alamo2005a} is based on the Mean Value Theorem and is referred to as the \emph{mean value extension}, while the approach in \cite{Combastel2005} uses a first-order Taylor expansion with a rigorous remainder bound, and is referred to as the \emph{first-order Taylor extension}. Updates are then achieved by methods for outer-approximating the intersection of a zonotope with a strip (i.e., a linear measurement with bounded error). An alternative zonotope-based prediction step using DC programming is proposed in \cite{Alamo2008}, but with the same update as in \cite{Alamo2005a}. Even for linear measurements, the symmetry of zonotopes is known to cause significant errors in the update step \cite{Scott2016}. General convex polytopes in halfspace representation are used in \cite{Wan2018} to enable an exact update. Prediction is then done by representing the polytope as an intersection of zonotopes and applying the mean value extension. Unfortunately, conversion between these representations is computationally demanding, and the increasing complexity of the zonotope bundle with time is not addressed. Constrained zonotopes have recently been applied to nonlinear state estimation in \cite{Rego2018} and shown to provide much higher accuracy than existing zonotopic methods. Nevertheless, the method in \cite{Rego2018} uses an interval partitioning scheme and is therefore intractable for high-dimensional systems. Finally, in an effort to overcome the limitations of convex sets, polynomial zonotopes are introduced in \cite{Althoff2013} and used for reachability analysis. However, update algorithms for polynomial zonotopes have not been developed.

In this context, the main contributions of this paper are two new methods for nonlinear set-valued state estimation based on constrained zonotopes. We follow the standard algorithmic steps typically used for set-valued state estimation (i.e., prediction, update, and reduction). For the prediction step, we use new generalizations of the mean value extension and first-order Taylor extension discussed above that enable constrained zonotopes, rather than zonotopes, to be effectively propagated through nonlinear discrete-time dynamics. Since this class of sets corresponds to an alternative representation of convex polytopes, these generalizations can be viewed as novel approaches for implicitly propagating convex polytopes through nonlinear mappings. The generalization of these methods to constrained zonotopes is not straightforward and requires significant modifications to the existing proofs. However, it results in much tighter prediction enclosures than existing zonotopic methods in some cases, as shown in the numerical examples. Moreover, using constrained zonotopes for prediction also enables the update step to be done much more effectively using the generalized intersection operation for constrained zonotopes. On the other hand, reduction becomes more complex than for zonotopes, leading to interesting trade-offs between accuracy and computational efficiency. We investigate these trade-offs both through numerical examples and by developing a detailed computational complexity analysis of the proposed methods and their zonotopic counterparts. To the best of our knowledge, such an analysis has not been conducted previously for either zonotopes or constrained zonotopes in the context of nonlinear set-valued estimation.

The remainder of the paper is organized as follows. The nonlinear set-valued state estimation problem is stated in Section \ref{sec:problem}. Essential mathematical background is presented in Section \ref{sec:preliminaries}, including a discussion of constrained zonotopes and their main properties. Section \ref{sec:estimation} develops the main results of the paper; namely, the proposed mean value and first-order Taylor extensions for constrained zonotopes, heuristics for selecting the point at which the approximation is performed, and the computational complexity analysis. Numerical examples are presented in Section \ref{sec:numericalexamples} to demonstrate the effectiveness of these extensions for set-valued state estimation of nonlinear discrete-time systems. Finally, Section \ref{sec:conclusion} concludes the paper.

\section{Problem formulation} \label{sec:problem}

Consider a class of discrete-time systems with nonlinear dynamics and linear measurements, described by
\begin{equation}
	\begin{aligned} \label{eq:system}
		\mathbf{x}_k & = \mathbf{f}(\mathbf{x}_{k-1}, \mathbf{u}_{k-1}, \mathbf{w}_{k-1}), \\
		\mathbf{y}_k & = \mathbf{C} \mathbf{x}_k + \mathbf{D}_u \mathbf{u}_k +\mathbf{D}_v \mathbf{v}_k,	
	\end{aligned}
\end{equation}
for $k \geq 1$, with $\mathbf{y}_0 = \mathbf{C} \mathbf{x}_0 + \mathbf{D}_u \mathbf{u}_0 +\mathbf{D}_v \mathbf{v}_0$, where $\mathbf{x}_k \in \realset^{n}$ denotes the system state, $\mathbf{u}_{k} \in \realset^{n_u}$ is a known input, $\mathbf{w}_k \in \realset^{n_w}$ is the process disturbance, $\mathbf{y}_k \in \realset^{n_y}$ is the measured output, and $\mathbf{v}_k \in \realset^{n_v}$ is the measurement uncertainty, with $\mbf{x}_0$ the initial state. The nonlinear mapping $\mathbf{f}$ is assumed to be of class $\mathcal{C}^2$, and the disturbances and uncertainties are assumed to be bounded, i.e., $\mathbf{w}_{k} \in W_k$ and $\mathbf{v}_k \in V_k$, where $W_k$ and $V_k$ are known compact sets.

This work proposes new methods to perform set-valued state estimation for nonlinear systems as in \eqref{eq:system}. The exact characterization of sets $X_k$ containing the evolution of the system states is very difficult in the nonlinear case, if not intractable \cite{Kieffer1998,Kuhn1998,Platzer2007}. Therefore, in the set-membership framework the objective is to enclose such sets as tightly as possible by guaranteed outer bounds $\hat{X}_k$ on the possible trajectories of the system states $\mathbf{x}_k$. Such outer bounds must be consistent with the previous estimate $\hat{X}_{k-1}$, known inputs $\mathbf{u}_{k-1}$, the current measurement $\mathbf{y}_k$, and also with the bounds on the disturbances and uncertainties $W_{k-1}$, $V_k$. Given an initial condition $\mathbf{x}_0 \in \hat{X}_0$, a common approach is to proceed through the well-known prediction-update algorithm, which consists in computing compact sets $\bar{X}_k$ and $\hat{X}_k$ such that
\begin{align}
	\bar{X}_k & \supseteq \{ \mathbf{f}(\mathbf{x}, \mathbf{u}_{k-1}, \mathbf{w}) : \mathbf{x} \in \hat{X}_{k-1}, \, \mathbf{w} \in W_{k-1} \}, \label{eq:prediction0}\\
	\hat{X}_k & \supseteq \{ \mathbf{x} \in \bar{X}_k : \mathbf{C} \mathbf{x} + \mathbf{D}_u \mathbf{u}_k + \mathbf{D}_v \mathbf{v} = \mathbf{y}_k, \, \mathbf{v} \in V_k \}, \label{eq:update0}
\end{align}
in which \eqref{eq:prediction0} is referred to as the \emph{prediction step}, and \eqref{eq:update0} as the \emph{update step}. 

Our goal is to obtain accurate outer bounds $\bar{X}_k$ and $\hat{X}_k$ according to \eqref{eq:prediction0} and \eqref{eq:update0}, respectively. Following these definitions, and considering the initial condition $\mathbf{x}_0 \in \hat{X}_0$, the property $\mathbf{x}_k \in \hat{X}_k$ is guaranteed by construction for all $k \geq \textred{1}$ \cite{Chisci1996,Le2013,Alamo2008}.

\section{Preliminaries} \label{sec:preliminaries}

\subsection{Set operations and constrained zonotopes}

A few common set operations are often used to compute enclosures satisfying \eqref{eq:prediction0} and \eqref{eq:update0} \cite{Le2013,Scott2016}. Consider $Z, W \subset \realset^{n}$, $\mathbf{R} \in \realset^{m \times n}$, and $Y \subset \realset^{m}$. Define the linear mapping, Minkowski sum, and generalized intersection, as
\begin{align}
	\mathbf{R}Z & \triangleq \{ \mathbf{R} \mathbf{z} : \mathbf{z} \in Z\}, \label{eq:limage}\\
	Z \oplus W & \triangleq \{ \mathbf{z} + \mathbf{w} : \mathbf{z} \in Z,\, \mathbf{w} \in W\}, \label{eq:msum}\\
	Z \cap_{\mathbf{R}} Y & \triangleq \{ \mathbf{z} \in Z : \mathbf{R} \mathbf{z} \in Y\}, \label{eq:intersection}
\end{align}
respectively. Using ellipsoids or parallelotopes, the linear mapping \eqref{eq:limage} can be computed exactly, but \eqref{eq:msum} and \eqref{eq:intersection} must be over-approximated \cite{Chisci1996,Schweppe1968}. For intervals, the Minkowski sum \eqref{eq:msum} is exact, but \eqref{eq:limage} and \eqref{eq:intersection} are conservative due to the wrapping effect\footnote{The generalized intersection in \eqref{eq:intersection} is not conservative when $\mathbf{R} = \eyenoarg$, which corresponds to the standard intersection $\cap$.}. In contrast, convex polytopes are closed under \eqref{eq:limage}--\eqref{eq:intersection}. Moreover, \eqref{eq:limage} and \eqref{eq:msum} can be computed efficiently in vertex representation (V-rep), and \eqref{eq:intersection} can be computed efficiently in half-space representation (H-rep). However, conversion between H-rep and V-rep is computationally expensive. Zonotopes \cite{Kuhn1998} allow \eqref{eq:limage} and \eqref{eq:msum} to be computed exactly and with low computational burden, but \eqref{eq:intersection} is not a zonotope in general and must be over-approximated \cite{Le2013,Bravo2006}.

Constrained zonotopes are an extension of zonotopes, recently proposed in \cite{Scott2016}, and are the class of sets of main interest in this work.

\begin{definition} \rm \label{def:czonotopes}
	A set $Z \subset \realset^n$ is a \emph{constrained zonotope} if there exists $(\mathbf{G}_z,\mathbf{c}_z,\mathbf{A}_z,\mathbf{b}_z) \in \realsetmat{n}{n_g} \times \realset^n \times \realsetmat{n_c}{n_g} \times \realset^{n_c}$ such that
	\begin{equation} \label{eq:czdefinition}
		Z = \left\{ \mathbf{c}_z + \mathbf{G}_z \bm{\xi} : \| \bm{\xi} \|_\infty \leq 1, \mathbf{A}_z \bm{\xi} = \mathbf{b}_z \right\}.
	\end{equation}	
\end{definition}

Equation \eqref{eq:czdefinition} is called \emph{constrained generator representation} (CG-rep), in which each column of $\mathbf{G}_z$ is a \emph{generator}, $\mathbf{c}_z$ is the \emph{center}, and $\mathbf{A}_z \bm{\xi} = \mathbf{b}_z$ are \emph{constraints}. In this work, we refer to $\bm{\xi}$ as the \emph{generator variables}. Let $B_\infty(\mathbf{A}_z,\mathbf{b}_z) = \{\bm{\xi} \in \realset^{n_g} : \ninf{\bm{\xi}} \leq 1,\,  \mathbf{A}_z \bm{\xi} = \mathbf{b}_z \}$ and $B_\infty^{n_g} = \{\bm{\xi} \in \realset^{n_g} : \ninf{\bm{\xi}} \leq 1\}$ be, respectively, the constrained unitary hypercube and the $n_g$-dimensional unitary hypercube\footnote{For simplicity, we drop the use of the superscript $n_g$ for $B_\infty(\mbf{A}_z,\mbf{b}_z)$. This dimension can be directly inferred from the number of columns of $\mbf{A}_z$.}. Then, a constrained zonotope $Z$ can alternatively be interpreted as an affine transformation of $B_\infty(\mathbf{A}_z,\mathbf{b}_z)$, given by $Z = \mathbf{c}_z \oplus \mathbf{G}_z B_\infty(\mathbf{A}_z,\mathbf{b}_z)$. Note that the linear equality constraints in \eqref{eq:czdefinition} allow constrained zonotopes to represent any convex polytope provided that the complexity of the CG-rep \eqref{eq:czdefinition} is not limited. In fact, $Z$ is a constrained zonotope iff it is a convex polytope \cite{Scott2016}. We use the compact notation $Z = \{\mathbf{G}_z, \mathbf{c}_z,\mathbf{A}_z,\mathbf{b}_z \}$ for constrained zonotopes, and $Z = \{\mathbf{G}_z, \mathbf{c}_z\}$ for zonotopes.

In addition to \eqref{eq:limage} and \eqref{eq:msum}, the intersection \eqref{eq:intersection} can also be computed exactly with constrained zonotopes. Let $Z = \{\mathbf{G}_z, \mathbf{c}_z, \mathbf{A}_z, \mathbf{b}_z\} \subset \realset^n$, $W = \{\mathbf{G}_w, \mathbf{c}_w, \mathbf{A}_w, \mathbf{b}_w\} \subset \realset^n$, $Y = \{\mathbf{G}_y, \mathbf{c}_y, \mathbf{A}_y, \mathbf{b}_y\} \subset \realset^m$, and $\mbf{R} \in \realsetmat{n}{m}$. The set operations \eqref{eq:limage}--\eqref{eq:intersection} are computed in CG-rep as
\begin{align}
	\mathbf{R}Z & = \left\{ \mathbf{R} \mathbf{G}_z, \mathbf{R} \mathbf{c}_z, \mathbf{A}_z, \mathbf{b}_z \right\}, \label{eq:czlimage}\\
	Z \oplus W & =\left\{ \begin{bmatrix} \mathbf{G}_z \,\; \mathbf{G}_w \end{bmatrix}, \mathbf{c}_z + \mathbf{c}_w, \begin{bmatrix} \mathbf{A}_z & \bm{0} \\ \bm{0} & \mathbf{A}_w \end{bmatrix}, \begin{bmatrix} \mathbf{b}_z \\ \mathbf{b}_w \end{bmatrix} \right\}\!, \label{eq:czmsum}\\
	Z \cap_{\mathbf{R}} Y & = \left\{ \begin{bmatrix} \mathbf{G}_z \,\; \bm{0} \end{bmatrix}, \mathbf{c}_z, \begin{bmatrix} \mathbf{A}_z & \bm{0} \\ \bm{0} & \mathbf{A}_y \\ \mathbf{R} \mathbf{G}_z & -\mathbf{G}_y \end{bmatrix}, \begin{bmatrix} \mathbf{b}_z \\ \mathbf{b}_y \\ \mathbf{c}_y - \mathbf{R} \mathbf{c}_z \end{bmatrix} \right\}. \label{eq:czintersection}
\end{align}

These operations can be performed efficiently and cause only a moderate increase in the complexity of the CG-rep. Other useful operations with constrained zonotopes are presented in the following. Property \ref{prope:ihull} provides a simple method for computing the interval hull of a constrained zonotope by solving $2n$ linear programs (LPs), while Proposition \ref{propo:closest} provides a way to obtain the closest point in a constrained zonotope to another point in space (in the 1-norm sense) through the solution of a single LP. For simplicity, the subscripts of the variables in \eqref{eq:czdefinition} will be suppressed henceforth when not necessary.

\begin{property}\rm  \cite{Scott2016,Rego2018} \label{prope:ihull}
    Let $Z = \{\mathbf{G}, \mathbf{c}, \mathbf{A}, \mathbf{b} \} \subset \realset^{n}$ and let $\mathbf{G}_j$ denote the $j$-th row of $\mathbf{G}$. The \emph{interval hull} $[\bm{\zeta}^\text{L},\,\bm{\zeta}^\text{U}] \supseteq Z$ is obtained by solving the following linear programs for each $j = 1,2,\dots,n$: 
	\begin{align*}
		\zeta_j^\text{L} & = \underset{\bm{\xi}}{\min} \left\{c_j + \mathbf{G}_j \bm{\xi} : \ninf{\bm{\xi}} \leq 1,\, \mathbf{A} \bm{\xi} = \mathbf{b} \right\}, \\
		\zeta_j^\text{U} & = \underset{\bm{\xi}}{\max} \left\{c_j + \mathbf{G}_j \bm{\xi} : \ninf{\bm{\xi}} \leq 1,\, \mathbf{A} \bm{\xi} = \mathbf{b} \right\}.
	\end{align*}
\end{property}

\begin{proposition} \rm \label{propo:closest} 
	Let $Z = \{\mathbf{G}, \mathbf{c}, \mathbf{A}, \mathbf{b} \} \subset \realset^{n}$ and $\mathbf{h} \in \realset^{n}$. A point $\mathbf{z} \in Z$ that minimizes $\none{\mathbf{z} - \mathbf{h}}$ is given by $\mathbf{z}^* = \mathbf{c} + \mathbf{G} \bm{\xi}^*$, where $\bm{\xi}^*$ is a solution to the linear program
	\begin{equation*}
	\underset{\bm{\xi}}{\min}~\| \mathbf{c} - \mathbf{h} + \mathbf{G} \bm{\xi} \|_1 \quad \text{\rm s.t.} \quad \mathbf{A}\bm{\xi} = \mathbf{b}, \quad \ninf{\bm{\xi}} \leq 1.
	\end{equation*}
\end{proposition}
\proof By definition,
\begin{align*}
\none{\mathbf{z}^* - \mathbf{h}} &= \none{\mathbf{c} - \mathbf{h} + \mathbf{G} \bm{\xi}^*} \\
&\leq \none{\mathbf{c} - \mathbf{h} + \mathbf{G} \bm{\xi}}, \quad \forall \bm{\xi}\in B_{\infty}(\mathbf{A},\mathbf{b}).
\end{align*}
But, for any $\mathbf{z}\in Z$, there exists $\bm{\xi}\in B_{\infty}(\mathbf{A},\mathbf{b})$ such that $\mathbf{z} = \mathbf{c} + \mathbf{G} \bm{\xi}$. Therefore, $\none{\mathbf{z}^* - \mathbf{h}} \leq \none{\mathbf{z} - \mathbf{h}}$, $\forall \mathbf{z}\in Z$.
\qed

The presence of the equality constraints in \eqref{eq:czdefinition} may result in $\mathbf{c} \notin Z$ (i.e., when $\bm{0} \notin B_\infty(\mathbf{A},\mathbf{b})$). Some of the techniques proposed in this paper require that $\mathbf{c} \in Z$ (Section \ref{sec:selectionofh}). To accommodate this, Proposition \ref{propo:center} provides a simple method for computing an alternative (more complex) CG-rep for $Z$ whose center is any desired point in space.

\begin{proposition}[Rescaling with desired center] \rm \label{propo:center}
	Let $Z = \{\mathbf{G},\mathbf{c},\mathbf{A},\mathbf{b}\} \subset \realset^n$ and  let $\tilde{\bm{\xi}}^\text{L}, \tilde{\bm{\xi}}^\text{U} \in \realset^{n_g}$ satisfy $\conball \subseteq [\tilde{\bm{\xi}}^\text{L},\,\tilde{\bm{\xi}}^\text{U}]$. Choose any desired center $\mathbf{h} \in \realset^n$ lying in the range of $\mathbf{G}$ and let $\xil, \xiu \in \realset^{n_g}$ be solutions to the linear program:
	\begin{align*}
    \underset{\xil,\xiu}{\min}~&\left\| \half (\xiu - \xil) \right\|_1 \\
	\text{s.t.} \quad & \mathbf{c} + \half \mathbf{G} (\xil + \xiu) = \mathbf{h}, \quad \xil \leq \tilde{\bm{\xi}}^\text{L}, \quad \xiu \geq \tilde{\bm{\xi}}^\text{U}.
	\end{align*}
    Letting $\xim = \half (\xil + \xiu)$ and $\mathbf{E}_r = \half\text{diag}(\xiu - \xil)$, an equivalent CG-rep of $Z$ with center $\mathbf{h}$ is given by
	\begin{align} 
	Z = \left\{\begin{bmatrix} \mathbf{G} \mathbf{E}_r & \bm{0} \end{bmatrix}, \mathbf{h}, \begin{bmatrix} \mathbf{A} \mathbf{E}_r & \bm{0} \\ \bm{0} & \mathbf{A} \\ \mathbf{G} \mathbf{E}_r & -\mathbf{G} \end{bmatrix} , \begin{bmatrix}\mathbf{b} - \mathbf{A} \xim  \\ \mathbf{b} \\ \mathbf{c} - \mathbf{h} \end{bmatrix}  \right\}. \label{eq:centerchange}
	\end{align}
\end{proposition}

\proof It is first shown that $Z$ is contained in the set 
\begin{equation*}
	\bar{Z} = \{ \mathbf{GE}_r,\,\mathbf{c} + \mathbf{G} \xim,\, \mathbf{AE}_r,\, \mathbf{b} - \mathbf{A} \xim \}.
	\end{equation*}
Choose any $\mathbf{z} \in Z$. There must exist $\bm{\xi} \in B_\infty (\mathbf{A},\mathbf{b})$ such that $\mathbf{z} = \mathbf{c} + \mathbf{G} \bm{\xi}$. Since $\bm{\xi} \in [\xil,\,\xiu]$, there must exist $\bm{\delta} \in B_\infty^{n_g}$ such that $ \bm{\xi} = \xim + \mathbf{E}_r\bm{\delta}$. Thus,
\begin{align*}
\mathbf{z} \in Z & \implies \exists \bm{\delta} \in B_\infty^{n_g} : \mathbf{z} = \mathbf{c} + \mathbf{G} (\xim + \mathbf{E}_r \bm{\delta}), \\ & \qquad \qquad \qquad \quad \; \mathbf{A}(\xim + \mathbf{E}_r \bm{\delta})= \mathbf{b} \implies \mathbf{z} \in \bar{Z}.
\end{align*}
Therefore, $Z \subseteq \bar{Z}$ and it is true that $Z=\bar{Z} \cap Z$. Since $\xil$ and $\xiu$ satisfy $\mathbf{c} + \half \mathbf{G} (\xil + \xiu) = \mathbf{c} + \mathbf{G} \xim = \mathbf{h}$, then representing $\bar{Z} \cap Z$ as in \eqref{eq:czintersection} gives \eqref{eq:centerchange}. \qed

\begin{remark} \rm \label{rem:rescalingtocenter}
    The linear program in Proposition \ref{propo:center} does not require $[\xil, \xiu] \subseteq B_\infty^{n_g}$. Therefore, the midpoint $\half (\xil + \xiu)$ can assume any desired value, and it is always possible to satisfy $\mathbf{c} + \half \mathbf{G} (\xil + \xiu) = \mathbf{h}$ if $\mathbf{h}$ is in the range of $\mathbf{G}$. Thus, the linear program is always feasible.
\end{remark}

\begin{remark} \rm The optimization problems in this section can be readily rewritten as standard form LPs by using additional decision variables and constraints \citep{Boyd2004}.
\end{remark}

\section{Nonlinear state estimation} \label{sec:estimation}

This section presents two new methods for set-valued state estimation of the class of nonlinear discrete-time systems described by \eqref{eq:system}. Focusing on the prediction step \eqref{eq:prediction0}, we address the problem of propagating a constrained zonotope $\hat{X}_{k-1}$ through a nonlinear mapping, with $\hat{X}_0$, $W_k$, and $V_k$ being constrained zonotopes as well. This is done by extending, in a consistent way, two existing approaches for propagating zonotopes through nonlinear mappings; namely, the \emph{mean value extension} in \cite{Alamo2005a,Kuhn1998} and the \emph{first-order Taylor extension} in \cite{Combastel2005}. The methods described in these works rely, respectively, on the Mean Value Theorem and Taylor's Theorem for the calculation of rigorous outer bounds of the range of the nonlinear mapping in order to obtain zonotopes enclosing the system trajectories. Both methods are based on intersection with strips for performing the update step \eqref{eq:update0}. The key advantage of our new extensions is that they allow the entire state estimation procedure to be done using constrained zonotopes in CG-rep. Therefore, the update \eqref{eq:update0} can be done by generalized intersection (with linear measurements), which is known to generate highly asymmetrical sets that cannot be accurately enclosed by ellipsoids, intervals, parallelotopes, or zonotopes. Using the methods developed below, such sets can be directly propagated to the next time step without prior simplification to a symmetric set. This overcomes a major source of conservatism in existing methods based on the aforementioned enclosures, while largely retaining the efficiency of computations with zonotopes. In addition, our methods expand the use of the important tools developed in \cite{Scott2016} to the class of nonlinear discrete-time systems described in Section \ref{sec:problem}. In the remainder of the paper, functions with set-valued arguments are consistently used to denote exact image of the set under the function; e.g., $\bm{\mu}(X,W)\triangleq\{\bm{\mu}(\mathbf{x},\mathbf{w}): \mathbf{x}\in X, \ \mathbf{w}\in W\}$.

The methods below make use of interval arithmetic in several places. Here we recall some of its main concepts. Let $\intvalset$ denote the set of compact intervals in $\realset$. Then $x \in \intvalset$ is a real compact set defined by $x = \{ a \in \realset : x^\lbound \leq a \leq x^\ubound \}$, with shorthand notation $[x^\lbound,x^\ubound]$. The midpoint and radius are defined by $\midpoint{x} = \half(x^\ubound + x^\lbound)$ and $\rad{x} = \half(x^\ubound - x^\lbound)$. The diameter is $\diam{x} = 2\rad{x}$. Interval vectors and matrices are defined by $\{ \mathbf{a} \in \realset^n : a_i^\lbound \leq a_i \leq a_i^\ubound\}$ and $\{ \mathbf{A} \in \realset^{n \times m}: A_{ij}^\lbound \leq A_{ij} \leq A_{ij}^\ubound\}$, respectively, with midpoint and radius defined component-wise. $\iextension{\mathbf{f}(X)}$ denotes an interval enclosure of a vector valued function $\mathbf{f}$ over $X \subset \realset^n$. The notation $\iextension{\mathbf{f}(X)}$ is used even when $X$ is not an interval. In this case, it is assumed that the interval hull of $X$ is employed in the operation. Inclusion functions satisfy $\mathbf{f}(X) \subseteq \iextension{\mathbf{f}(X)}$. See \cite{Moore2009} for definitions of basic interval arithmetic operations and examples.

In addition, the following notations are defined to be used in our proofs. Let $\bm{\kappa}$ be a function of class $\mathcal{C}^2$ and $\mbf{z}$ denote its argument. Then, $\kappa_q$ denotes the $q$-th component of $\bm{\kappa}$, $\nabla \bm{\kappa}$ denotes the gradient of $\bm{\kappa}$, and $\mathbf{H} \kappa_q$ is an upper triangular matrix describing half of the Hessian of $\kappa_q$. Specifically, $H_{ii} \kappa_q = (1/2) \partial^2 \kappa_q/\partial z_i^2$, $H_{ij} \kappa_q = \partial^2 \kappa_q/\partial z_i \partial z_j$ for $i<j$, and $H_{ij} \kappa_q = 0$ for $i>j$.

\subsection{Mean value extension}

This section presents the first of two new methods for enclosing the range of a nonlinear function $\bm{\mu}$ over a set of inputs described by constrained zonotopes. This method is referred to as the \emph{mean value extension} of $\bm{\mu}$ (because it relies on the Mean Value Theorem), and is a consistent generalization of the method for zonotopes in \cite{Alamo2005a}. Due to significant differences with respect to its zonotopic counterpart, a new theorem (Theorem \ref{the:pred2}) together with a detailed proof is provided for the new method.

The method in \cite{Alamo2005a} relies on a \emph{zonotope inclusion} operator that computes a zonotopic enclosure of the product of an interval matrix with a unitary box. We first generalize this operator to constrained zonotopes.

\begin{theorem} \rm \label{the:gzielim}
Let $X = \mathbf{p} \oplus \mathbf{M} B_\infty (\mathbf{A},\mathbf{b}) \subset \realset^{m}$ be a constrained zonotope with $n_g$ generators and $n_c$ constraints, let $\mathbf{J} \in \intvalsetmat{n}{m}$ be an interval matrix, and consider the set $S = \mathbf{J} X \triangleq \{\hat{\mathbf{J}} \mathbf{x} : \hat{\mathbf{J}} \in \mathbf{J}, \mathbf{x} \in X\} ~\subset \realset^n$.
Let $\bar{X} = \bar{\mathbf{p}} \oplus \bar{\mathbf{M}} B_\infty^{\bar{n}_g}$ be a zonotope satisfying $X \subseteq \bar{X}$, let $\mathbf{m}$ be an interval vector such that $\mathbf{m} \supseteq (\mathbf{J} - \midpoint{\mathbf{J}}) \bar{\mathbf{p}}$ and $\midpoint{\mathbf{m}} = \bm{0}$, and let $\mathbf{P} \in \realsetmat{n}{n}$ be a diagonal matrix defined by
	\begin{equation} \label{eq:czinclusionP}
	P_{ii} = \half \diam{m_i} +  \half \sum_{j=1}^{\bar{n}_g} \sum_{k=1}^{m} \diam{J_{ik}} |\bar{M}_{kj}|,
	\end{equation}
	for all $i=1,2,\dots,n$. Then $S$ is contained in the \emph{CZ-inclusion}
	\begin{equation} \label{eq:czinclusion}
		S \subseteq \gzinclusion(\mathbf{J},X) \triangleq \midpoint{\mathbf{J}}X \oplus \mathbf{P}B_\infty^n.
	\end{equation}
\end{theorem}

\proof 
Choose any $\mathbf{s} \in S$. It will be shown that $\mathbf{s} \in \gzinclusion(\mathbf{J},X)$. By the definition of $S$, there must exist $\mathbf{x} \in X$ and $\hat{\mathbf{J}} \in \mathbf{J}$ such that $\mathbf{s} = \hat{\mathbf{J}} \mathbf{x}$. Adding and subtracting $\midpoint{\mathbf{J}} \mathbf{x}$,
\begin{equation*}
	\mathbf{s} = \midpoint{\mathbf{J}} \mathbf{x} + (\hat{\mathbf{J}} - \midpoint{\mathbf{J}}) \mathbf{x}.	
\end{equation*}
Since $\mathbf{x}\in X\subseteq \bar{X}$, there exists $\bm{\gamma} \in B_\infty^{\bar{n}_g}$ such that $\mathbf{x} = \bar{\mathbf{p}} + \bar{\mathbf{M}} \bm{\gamma}$. Therefore, $\mathbf{s} = \midpoint{\mathbf{J}}\mathbf{x} + (\hat{\mathbf{J}} - \midpoint{\mathbf{J}}) (\bar{\mathbf{p}} + \bar{\mathbf{M}} \bm{\gamma}).$ By the choice of $\mathbf{m}$, there must exist $\hat{\mathbf{m}}\in\mathbf{m}$ such that
\begin{align}
\label{Eq: CZ Inc Proof - z equality}
\mathbf{s} &= \midpoint{\mathbf{J}}\mathbf{x} + \hat{\mathbf{m}} + (\hat{\mathbf{J}} - \midpoint{\mathbf{J}}) \bar{\mathbf{M}} \bm{\gamma}.
\end{align}
Let $\bm{\eta}= \hat{\mathbf{m}}+(\hat{\mathbf{J}} - \midpoint{\mathbf{J}}) \bar{\mathbf{M}} \bm{\gamma}$. Then 
\begin{align*}
\eta_{i} &= \hat{m}_i+\sum_{j=1}^{\bar{n}_g} ((\hat{\mathbf{J}} - \midpoint{\mathbf{J}}) \bar{\mathbf{M}})_{ij} \gamma_j, \\
&= \hat{m}_i+\sum_{j=1}^{\bar{n}_g} \left(\sum_{k=1}^{m} (\hat{J}_{ik}-\midpoint{J_{ik}}) \bar{M}_{kj}\right) \gamma_j.
\end{align*}
By the triangle inequality and the fact that $|\gamma_j|\leq 1$,
\begin{align*}
|\eta_{i}| &\leq |\hat{m}_i|+\sum_{j=1}^{\bar{n}_g} \left(\sum_{k=1}^{m} |(\hat{J}_{ik}-\midpoint{J_{ik}})|| \bar{M}_{kj}|\right)|\gamma_j|, \\
&\leq \half \diam{m_i} +  \half \sum_{j=1}^{\bar{n}_g} \sum_{k=1}^{m} \diam{J_{ik}} |\bar{M}_{kj}|.
\end{align*}
Therefore, $\bm{\eta}\in \mathbf{P}B_\infty^n$. From \eqref{Eq: CZ Inc Proof - z equality}, this implies that
\begin{equation*}
\mathbf{s} = \midpoint{\mathbf{J}}\mathbf{x} + \bm{\eta}
\in \midpoint{\mathbf{J}} X \oplus \mathbf{P} B_\infty^n = \gzinclusion(\mathbf{J},X).
\end{equation*}
Thus $S \subseteq \gzinclusion(\mathbf{J},X)$.
\qed

\begin{remark} \rm  \label{rem:contrelim}
	In Theorem \ref{the:gzielim}, a zonotope satisfying $X \subseteq \bar{X}$ can be easily obtained by performing $n_c$ iterated constraint eliminations on $X$ using the method in \cite{Scott2016}. Moreover, $\mathbf{m}$ can be obtained by simply evaluating $(\mathbf{J} - \midpoint{\mathbf{J}}) \bar{\mathbf{p}}$ with interval arithmetic. These methods are used in this work. Finally, the enclosure \eqref{eq:czinclusion} has $n_g + n$ generators and $n_c$ constraints.
\end{remark}

The following theorem provides the mean value extension for constrained zonotopes. 

\begin{theorem} \rm \label{the:pred2}
     Let $\bm{\mu} : \realset^n \times \realset^{n_w} \to \realset^n$ be continuously differentiable and $\nabla_x \bm{\mu}$ denote the gradient of $\bm{\mu}$ with respect to its first argument. Let $X\subset \realset^n$ and $W \subset \realset^{n_w}$ be constrained zonotopes and choose any $\mathbf{h} \in X$. If $Z$ is a constrained zonotope such that $\bm{\mu}(\mathbf{h},W) \subseteq Z$ and $\mathbf{J} \in \intvalsetmat{n}{n}$ is an interval matrix satisfying $\nabla^T_x \bm{\mu}(X, W)\subseteq \mathbf{J}$, then $\bm{\mu}(X,W) \subseteq Z \oplus \gzinclusion\left(\mathbf{J},  X - \mathbf{h} \right)$.
\end{theorem}

\proof Choose any $(\mathbf{x},\mathbf{w})\in X\times W$. It will be shown that $\bm{\mu}(\mathbf{x},\mathbf{w}) \in Z \oplus \gzinclusion\left(\mathbf{J},  X - \mathbf{h} \right)$. For any $i=1,2,\ldots,n$, the Mean Value Theorem ensures that $\exists \bm{\delta}^{[i]}\in X$ such that
\begin{align*}
	\mu_i(\mathbf{x},\mathbf{w}) &= \mu_i(\mathbf{h},\mathbf{w}) + \nabla^T_x \mu_i(\bm{\delta}^{[i]},\mathbf{w}) (\mathbf{x} - \mathbf{h}).
\end{align*}
But the vector $\nabla^T_x \mu_i(\bm{\delta}^{[i]},\mathbf{w})$ is contained in the $i$-th row of $\mathbf{J}$ by hypothesis, and since this is true for all $i=1,2,\ldots,n$, there exists a real matrix $\hat{\mathbf{J}}\in \mathbf{J}$ such that $\bm{\mu}(\mathbf{x},\mathbf{w}) = \bm{\mu}(\mathbf{h},\mathbf{w}) + \hat{\mathbf{J}} (\mathbf{x} - \mathbf{h}).$ 

By Theorem \ref{the:gzielim} and the choice of $Z$, it follows that $\bm{\mu}(\mathbf{x},\mathbf{w}) \in Z \oplus \gzinclusion\left(\mathbf{J}, X - \mathbf{h}\right)$, as desired. \qed

\begin{remark} \rm \label{rem:affine}
The interval matrix $\mathbf{J}$ required by Theorem \ref{the:pred2} can be obtained by computing the interval hulls of $X$ and $W$ as in Property \ref{prope:ihull} and then bounding $\nabla^T_x \bm{\mu}(X, W)$ using interval arithmetic. Similarly, the constrained zonotope $Z \supseteq \bm{\mu}(\mathbf{h},W)$ can be obtained by bounding $\bm{\mu}(\mathbf{h},W)$ with interval arithmetic. Alternatively, another mean value extension can be applied around some $\mathbf{h}_w \in W$ to obtain $\bm{\mu}(\mathbf{h},W) \subseteq Z\triangleq \bm{\mu}(\mathbf{h},\mathbf{h}_w) \oplus \gzinclusion\left(\mathbf{J}_w,W - \mathbf{h}_w\right)$, where $\mathbf{J}_w$ is an interval enclosure of $\nabla^T_w\bm{\mu}(\mathbf{h}, W)$. Finally, if $\bm{\mu}$ is affine in $\mathbf{w}$, i.e, $\bm{\mu}(\mathbf{x},\mathbf{w}) \triangleq \bm{\beta}_x(\mathbf{x}) + \mbf{B}_w(\mathbf{x}) \mathbf{w}$, then an exact enclosure of $\bm{\mu}(\mathbf{h},W)$ is $Z = \bm{\beta}_x(\mathbf{h}) \oplus \mbf{B}_w(\mathbf{h}) W$.
\end{remark}

Since the CG-rep \eqref{eq:czdefinition} is an alternative representation for convex polytopes \cite{Scott2016}, the mean value extension developed in Theorem \ref{the:pred2} provides a new method for propagating convex polytopes implicitly through nonlinear mappings. A related approach can be found in \cite{Wan2018}, where convex polytopes are represented by intersections of zonotopes (i.e., zonotope bundles \cite{Althoff2011}). However, while effective complexity reduction algorithms are available for constrained zonotopes \cite{Scott2016}, efficient methods for complexity control of zonotope bundles have not yet been proposed. 

\begin{remark} \rm \label{rem:meanvaluecomplexity}
The enclosure obtained in Theorem \ref{the:pred2} has at most $n_g + n_{g_w} + 2n$ generators and $n_c + n_{c_w}$ constraints (considering $Z$ computed as in the alternatives presented in Remark \ref{rem:affine}), with $n_g$ and $n_{g_w}$ denoting the number of generators of $X$ and $W$, and $n_c$ and $n_{c_w}$ the number of constraints, respectively. Thus, the complexity of the resulting set grows linearly with respect to the number of constraints and generators.
\end{remark}

\subsection{First-order Taylor extension}
 
This section presents the second new method for enclosing the range of a nonlinear function $\bm{\eta}$ over a set of inputs described by constrained zonotopes. This method is referred to as the \emph{first-order Taylor extension} of $\bm{\eta}$ (because it relies on a first-order Taylor expansion with a rigorous remainder bound), and is a consistent generalization of the method for zonotopes in \cite{Combastel2005}. In contrast to Theorem \ref{the:pred2}, for the sake of simplicity of the proof, this function has only one argument. Even so, it is possible to consider both states and process uncertainties by concatenating them into a single vector. Due to substantial changes with respect to the zonotopic method, the new approach comes with a new theorem (Theorem \ref{thm:firstorder}) and a detailed proof. In the main result below, $(\cdot)_{i,:}$ denotes the $i$-th row of a matrix, and $(\cdot)_{ij}$ denotes the element from its $i$-th row and $j$-th column.

\begin{theorem} \label{thm:firstorder} \rm
    Let $\bm{\eta}: \realset^{m} \to \realset^{n}$ be of class $\mathcal{C}^2$ and $\mathbf{z} \in \realset^{m}$ denote its argument. Let $Z = \{\mathbf{G}, \mathbf{c}, \mathbf{A}, \mathbf{b}\} \subset \realset^{m}$ be a constrained zonotope with $m_g$ generators and $m_c$ constraints. For each $q = 1,2,\dots,n$, let $\mathbf{Q}^{[q]}\in\mathbb{IR}^{m\times m}$ and $\tilde{\mathbf{Q}}^{[q]}\in\mathbb{IR}^{m_g\times m_g}$ be interval matrices satisfying $\mathbf{Q}^{[q]} \supseteq \mathbf{H} \eta_q (Z)$ and $\tilde{\mathbf{Q}}^{[q]} \supseteq \mathbf{G}^T \mathbf{Q}^{[q]} \mathbf{G}$. Moreover, define        
	\begin{align*}
		& \tilde{c}_q = \trace{\midpoint{\tilde{\mathbf{Q}}^{[q]}}}/2, \\
		& \tilde{\mathbf{G}}_{q,:} = \big[ \cdots \,\; \underbrace{\midpoint{\tilde{Q}^{[q]}_{ii}}/2}_{\forall i} \,\; \cdots \,\; \underbrace{\left(\midpoint{\tilde{Q}^{[q]}_{ij}} + \midpoint{\tilde{Q}^{[q]}_{ji}}\right)}_{\forall i<j} \,\; \cdots \big],\\
		& \tilde{\mathbf{G}}_{\mathbf{d}} = \text{diag}(\mathbf{d}), \quad d_q = \sum_{i,j} \left| \rad{\tilde{Q}^{[q]}_{ij}} \right|, \\
		& \tilde{\mathbf{A}} = \left[ \tilde{\mathbf{A}}_{\bm{\zeta}} \,\; \tilde{\mathbf{A}}_{\bm{\xi}} \,\; \zeros{\frac{m_c}{2}(1+m_c)}{n} \right],\\
		& \tilde{\mathbf{A}}_{\bm{\zeta}} = \arrowmatrix{\begin{matrix} & \vdots & \\ \cdots & \half A_{ri} A_{si} & \cdots \\ & \vdots & \end{matrix}}{\forall i}{\forall r \leq s}, \\
		& \tilde{\mathbf{A}}_{\bm{\xi}} = 
		\arrowmatrix{\begin{matrix} & \vdots & \\ \cdots & A_{ri} A_{sj} + A_{rj} A_{si} & \cdots \\ & \vdots & \end{matrix}}{\forall i < j}{\forall r \leq s}, \\
		& \tilde{\mathbf{b}} = \downarrowmatrix{\begin{matrix} \vdots \\ b_{r} b_{s} - \half \sum_i A_{ri} A_{si} \\ \vdots\end{matrix}}{\forall r \leq s},
	\end{align*}
    with indices $i,j = 1,2,\dots,m_g$ and $r,s = 1,2,\dots,m_c$. Finally, choose any $\mathbf{h} \in Z$ and let $\mathbf{L}\in\mathbb{IR}^{n\times m}$ be an interval matrix satisfying $\mathbf{L}_{q,:} \supseteq (\mathbf{c} - \mathbf{h})^T \mathbf{Q}^{[q]}$ for all $q = 1,\dots,n$. Then,
	\begin{equation} \label{eq:firstorderextension}
		\bm{\eta}(Z) \subseteq \bm{\eta}(\mathbf{h}) \oplus \nabla^T \bm{\eta}(\mathbf{h})(Z - \mathbf{h}) \oplus R,
	\end{equation}
    where $R = \tilde{\mathbf{c}} \oplus \left[ \tilde{\mathbf{G}} \,\; \tilde{\mathbf{G}}_{\mathbf{d}} \right] B_\infty(\tilde{\mathbf{A}}, \tilde{\mathbf{b}}) \oplus \gzinclusion (\mathbf{L}, (\mathbf{c} - \mathbf{h}) \oplus 2\mathbf{G} B_\infty(\mathbf{A},\mathbf{b}) )$.
\end{theorem}

\proof
Choose any $\mathbf{z}\in Z$ and $q\in\{1,\ldots,n\}$. By Taylor's theorem applied to $\eta_q$ with reference point $\mathbf{h}$, there must exist $\bm{\Gamma}^{[q]} \in \mathbf{H} \eta_q (Z) \subseteq \mathbf{Q}^{[q]}$ such that \footnote{Let $\bm{\Upsilon}^{[q]}$ belong to the standard Hessian matrix of $\eta_q (Z)$. Then, $(1/2)(\mathbf{z} - \mathbf{h})^T \bm{\Upsilon}^{[q]} (\mathbf{z} - \mathbf{h}) =  (\mathbf{z} - \mathbf{h})^T \bm{\Gamma}^{[q]} (\mathbf{z} - \mathbf{h})$ holds. See \citep{Combastel2005} for a motivation on this approach.}
\begin{equation*}
	\eta_q(\mathbf{z}) = \eta_q (\mathbf{h}) + \nabla^T \eta_q (\mathbf{h})(\mathbf{z} - \mathbf{h}) + (\mathbf{z} - \mathbf{h})^T \bm{\Gamma}^{[q]} (\mathbf{z} - \mathbf{h}).
\end{equation*}
Since $\mathbf{z}\in Z$, there must exist $\bm{\xi} \in B_\infty(\mathbf{A},\mathbf{b})$ such that $\mathbf{z}=\mathbf{c} + \mathbf{G}\bm{\xi}$. Thus, defining $\mathbf{p} = \mathbf{c} - \mathbf{h}$ for brevity,
\begin{equation*}
	\eta_q(\mathbf{z}) = \eta_q (\mathbf{h}) + \nabla^T \eta_q (\mathbf{h})(\mathbf{p} + \mathbf{G}\bm{\xi}) + (\mathbf{p} + \mathbf{G}\bm{\xi})^T \bm{\Gamma}^{[q]}  (\mathbf{p} + \mathbf{G}\bm{\xi}).
\end{equation*}
Expanding the product $(\mathbf{p} + \mathbf{G}\bm{\xi})^T \bm{\Gamma}^{[q]}  (\mathbf{p} + \mathbf{G}\bm{\xi})$ yields $\mathbf{p}^T \bm{\Gamma}^{[q]} (\mathbf{p} + 2 \mathbf{G} \bm{\xi}) + \bm{\xi}^T \tilde{\bm{\Gamma}}^{[q]} \bm{\xi}$, 
with $\tilde{\bm{\Gamma}}^{[q]} = \mathbf{G}^T \bm{\Gamma}^{[q]} \mathbf{G} \in \tilde{\mathbf{Q}}^{[q]}$.
Since $\tilde{\bm{\Gamma}}^{[q]} \in \tilde{\mathbf{Q}}^{[q]}$, it follows that $\tilde{\Gamma}^{[q]}_{ij} = \midpoint{\tilde{Q}^{[q]}_{ij}} + \rad{\tilde{Q}^{[q]}_{ij}} \Lambda^{[q]}_{ij}$ for some $\Lambda^{[q]}_{ij} \in B_\infty^1$. Additionally, $\xi_i \in [-1,1]$ implies that $\xi_i^2 \in [0, 1]$, and hence $\xi_i^2 = \half + \half \zeta_i$ for some $\zeta_i \in [-1,1]$. Considering these two facts,
\begin{align*}
	\bm{\xi}^T \tilde{\bm{\Gamma}}^{[q]} \bm{\xi} & = \half \sum_i \midpoint{\tilde{Q}^{[q]}_{ii}} + \half \sum_i \midpoint{\tilde{Q}^{[q]}_{ii}} \zeta_i \\
	& \quad + \sum_{i<j}(\midpoint{\tilde{Q}^{[q]}_{ij}} + \midpoint{\tilde{Q}^{[q]}_{ji}}) \xi_i \xi_j \\ & \quad + \sum_{i,j} \rad{\tilde{Q}^{[q]}_{ij}} \xi_i \xi_j \Lambda^{[q]}_{ij},
\end{align*}
where the third summation results from the fact that $\xi_i \xi_j = \xi_j \xi_i$. Thus, by defining the new generator variables
\begin{equation*} 
	\bar{\bm{\xi}} = \big[ \,\; \cdots \,\; \underbrace{\zeta_i}_{\forall i} \,\; \cdots \,\; \underbrace{\xi_i \xi_j}_{\forall i<j} \,\; \cdots \,\; \underbrace{\xi_i \xi_j \Lambda^{[q]}_{ij}}_{\forall i,j,q} \,\; \cdots \,\; \big]^T,
\end{equation*}
with $i,j = 1,2,\dots,m_g$, $q = 1,2,\dots,n$, we have that $\bm{\xi}^T \tilde{\bm{\Gamma}}^{[q]} \bm{\xi} = \tilde{c}_q + [\tilde{\mathbf{G}} \,\; \bar{\mathbf{G}}_{\mathbf{d}}]_{q,:} \bar{\bm{\xi}}$, where
$\bar{\mathbf{G}}_{\mathbf{d}} = \text{blkdiag}(\mathbf{N}^{[1]}, \mathbf{N}^{[2]}, \dots, \mathbf{N}^{[n]})$\footnote{In this work $\text{blkdiag}(\mathbf{A},\mathbf{B},\dots)$ denotes a block diagonal matrix with blocks $\mathbf{A},\mathbf{B},\dots$.}, 
\begin{equation*}
\mathbf{N}^{[q]} = \big[ \,\; \cdots \,\; \underbrace{\rad{\tilde{Q}^{[q]}_{ij}}}_{\forall i,j} \,\; \cdots \,\; \big] \in \realsetmat{1}{m_g^2}.
\end{equation*}

Therefore, we have established that $\eta_q(\mathbf{z}) = \eta_q (\mathbf{h}) + \nabla^T \eta_q (\mathbf{h})(\mathbf{z} - \mathbf{h}) + \mathbf{p}^T \bm{\Gamma}^{[q]} (\mathbf{p} + 2 \mathbf{G} \bm{\xi}) + \tilde{c}_q + [\tilde{\mathbf{G}} \,\; \bar{\mathbf{G}}_{\mathbf{d}}]_{q,:} \bar{\bm{\xi}}.$ This holds for every $q = 1,2,\dots,n$. Moreover, $\mathbf{L}$ satisfies $\mathbf{L}_{q,:} \supseteq \mathbf{p}^T \mathbf{Q}^{[q]}$ for all $q = 1,2,\dots,n$ by definition, so there must exist $\hat{\mathbf{L}}\in\mathbf{L}$ such that $\hat{\mathbf{L}}_{q,:}=\mathbf{p}^T \bm{\Gamma}^{[q]}$ for all $q = 1,2,\dots,n$. Therefore,
\begin{equation}
\label{Eq: eta vector expansion}
	\bm{\eta}(\mathbf{z}) = \bm{\eta}(\mathbf{h}) + \nabla^T \bm{\eta} (\mathbf{h})(\mathbf{z} - \mathbf{h})+ \hat{\mathbf{L}}(\mathbf{p} + 2 \mathbf{G} \bm{\xi}) + \tilde{\mathbf{c}} + [\tilde{\mathbf{G}} \,\; \bar{\mathbf{G}}_{\mathbf{d}}] \bar{\bm{\xi}}.
\end{equation}

Furthermore, the equality constraints $\mathbf{A} \bm{\xi} = \mathbf{b}$ imply that $\mathbf{A} \bm{\xi} \bm{\xi}^T \mathbf{A}^T = \mathbf{b} \mathbf{b}^T$. 
Thus, considering $\xi_i^2 = \half + \half \zeta_i$, the $r$-th row and $s$-th column of this matrix equality yields
\begin{equation*}
	\half \sum_i A_{ri} A_{si} \zeta_i + \sum_{i<j} (A_{ri} A_{sj} + A_{rj} A_{si}) \xi_i \xi_j = b_r b_s - \half \sum_i A_{ri} A_{si},
\end{equation*}
with $r,s = 1,2,\dots,m_c$.
Such constraints are linear in $\bar{\bm{\xi}}$, and non-repeating for $r \leq s$, therefore $\bar{\mathbf{A}} \bar{\bm{\xi}} = \tilde{\mathbf{b}}$ holds, where $\bar{\mathbf{A}} = [\tilde{\mathbf{A}}_{\bm{\zeta}} \,\; \tilde{\mathbf{A}}_{\bm{\xi}} \,\; \bm{0}_{\tilde{m}_c \times n m_g^2}]$, with $\tilde{m}_c = \frac{m_c}{2}(1+m_c)$. Hence $\bm{\xi} \in B_\infty(\mathbf{A},\mathbf{b}) \implies \bar{\bm{\xi}} \in B_\infty(\bar{\mathbf{A}},\tilde{\mathbf{b}})$. %
Combining this with \eqref{Eq: eta vector expansion}, we have proven the enclosure $\bm{\eta}(Z) \subseteq \bm{\eta} (\mathbf{h}) \oplus \nabla^T \bm{\eta} (\mathbf{h})(Z - \mathbf{h}) \oplus \mathbf{L} (\mathbf{p} \oplus 2 \mathbf{G} B_\infty(\mathbf{A},\mathbf{b})) \oplus \tilde{\mathbf{c}} \oplus [\tilde{\mathbf{G}} \,\; \bar{\mathbf{G}}_{\mathbf{d}}] B_\infty (\bar{\mathbf{A}}, \tilde{\mathbf{b}})$.

In fact, this enclosure can be greatly simplified by noting that the columns of $\bar{\mathbf{A}}$ corresponding to the variables $[\,\; \cdots \,\; \xi_i \xi_j \Lambda^{[q]}_{ij}\,\; \cdots \,\;]$ are all zero, and hence 
\begin{equation*}
[\tilde{\mathbf{G}} \,\; \bar{\mathbf{G}}_{\mathbf{d}}] B_\infty (\bar{\mathbf{A}}, \tilde{\mathbf{b}}) = \tilde{\mathbf{G}} B_\infty ([\tilde{\mathbf{A}}_{\bm{\zeta}} \,\; \tilde{\mathbf{A}}_{\bm{\xi}}], \tilde{\mathbf{b}}) \oplus \bar{\mathbf{G}}_{\mathbf{d}} B_\infty^{nm_g^2}.
\end{equation*}
Since $\bar{\mathbf{G}}_{\mathbf{d}}$ is block diagonal and each $\mathbf{N}^{[q]}$ is a row vector, $\bar{\mathbf{G}}_{\mathbf{d}} B_\infty^{nm_g^2}$ is an interval, and is equivalent to $\tilde{\mathbf{G}}_{\mathbf{d}} B_\infty^{n}$, with $\tilde{\mathbf{G}}_{\mathbf{d}}$ defined as in the statement of the theorem. Thus,
\begin{align*}
[\tilde{\mathbf{G}} \,\; \bar{\mathbf{G}}_{\mathbf{d}}] B_\infty (\bar{\mathbf{A}}, \tilde{\mathbf{b}}) &= \tilde{\mathbf{G}} B_\infty ([\tilde{\mathbf{A}}_{\bm{\zeta}} \,\; \tilde{\mathbf{A}}_{\bm{\xi}}], \tilde{\mathbf{b}}) \oplus \tilde{\mathbf{G}}_{\mathbf{d}} B_\infty^{n}, \\
&= [\tilde{\mathbf{G}} \,\; \tilde{\mathbf{G}}_{\mathbf{d}}] B_\infty ([\tilde{\mathbf{A}}_{\bm{\zeta}} \,\; \tilde{\mathbf{A}}_{\bm{\xi}} \,\; \mathbf{0}_{\tilde{m}_c\times n}], \tilde{\mathbf{b}}) \\ & = [\tilde{\mathbf{G}} \,\; \tilde{\mathbf{G}}_{\mathbf{d}}] B_\infty (\tilde{\mathbf{A}}, \tilde{\mathbf{b}}).
\end{align*}
Therefore, $\bm{\eta}(Z) \subseteq \bm{\eta} (\mathbf{h}) \oplus \nabla^T \bm{\eta} (\mathbf{h})(Z - \mathbf{h})
    \oplus \mathbf{L} (\mathbf{p} \oplus 2 \mathbf{G} B_\infty(\mathbf{A},\mathbf{b})) \oplus \tilde{\mathbf{c}} \oplus [\tilde{\mathbf{G}} \,\; \tilde{\mathbf{G}}_{\mathbf{d}}] B_\infty (\tilde{\mathbf{A}}, \tilde{\mathbf{b}}),$ 
and \eqref{eq:firstorderextension} follows immediately from the definition of $R$. \qed

\begin{remark} \rm
	Regarding the definitions of $\tilde{\mathbf{G}}$, $\tilde{\mathbf{A}}_{\bm{\zeta}}$, $\tilde{\mathbf{A}}_{\bm{\xi}}$, and $\tilde{\mathbf{b}}$ in Theorem \ref{thm:firstorder}, the ordering of the indices $i<j$ and $r\leq s$ is irrelevant, as long as it is the same for all variables.
\end{remark}
	
\begin{remark} \rm
	The interval matrices $\mathbf{Q}^{[q]}$ required by Theorem \ref{thm:firstorder} can be obtained by computing the interval hull of $Z$ (Property \ref{prope:ihull}) and then bounding $\mathbf{H} \eta_q (Z)$ using interval arithmetic. Moreover, $\tilde{\mathbf{Q}}^{[q]}$ and $\mathbf{L}$ can be obtained by evaluating $\mathbf{G}^T \mathbf{Q}^{[q]} \mathbf{G}$ and $(\mathbf{c} - \mathbf{h})^T \mathbf{Q}^{[q]}$ using interval arithmetic.
\end{remark}

As stated before, process disturbances can be taken into account in \eqref{eq:firstorderextension} by considering the augmented vector $\mathbf{z} = (\mathbf{x},\mathbf{w})$ with $Z = X \times W \subset \realset^{n+n_w}$ and $\mathbf{h} = (\mathbf{h}_x,\mathbf{h}_w) \in Z$. With $X = \{\mathbf{G}_x, \mathbf{c}_x, \mathbf{A}_x, \mathbf{b}_x\}$ and $W = \{\mathbf{G}_w, \mathbf{c}_w, \mathbf{A}_w, \mathbf{b}_w\}$, the Cartesian product $Z$ is easily computed by
\begin{equation*}
	X \times W = \left\{ \begin{bmatrix} \mathbf{G}_x & \bm{0} \\ \bm{0} & \mathbf{G}_w \end{bmatrix}, \begin{bmatrix} \mathbf{c}_x \\ \mathbf{c}_w \end{bmatrix}, \begin{bmatrix} \mathbf{A}_x & \bm{0} \\ \bm{0} & \mathbf{A}_w \end{bmatrix}, \begin{bmatrix} \mathbf{b}_x \\ \mathbf{b}_w \end{bmatrix} \right\}.
\end{equation*}

\begin{remark} \rm \label{rem:firstordercomplexity}
In Theorem \ref{thm:firstorder}, $\tilde{\mbf{G}}$ has $\sum_{j=1}^{m_g}j = \frac{1}{2}m_g(m_g+1)$ columns, $\tilde{\mbf{G}}_{\mbf{d}} \in \realsetmat{n}{n}$, $\tilde{\mbf{A}}$ has $\sum_{s=1}^{m_c}s = \frac{1}{2}m_c(m_c+1)$ rows, and $\gzinclusion (\mathbf{L}, (\mathbf{c} - \mathbf{h}) \oplus 2\mathbf{G} B_\infty(\mathbf{A},\mathbf{b}) )$ has $m_g+n$ generators and $m_c$ constraints (Remark \ref{rem:contrelim}). Therefore, the resulting enclosure in \eqref{eq:firstorderextension} has $\frac{1}{2}m_g^2 + \frac{5}{2}m_g + 2n$ generators and $\frac{1}{2}m_c^2 + \frac{5}{2}m_c$ constraints. If $Z = X \times W$, then the enclosure has $\frac{1}{2}(n_g + n_{g_w})^2 + \frac{5}{2}(n_g + n_{g_w}) + 2n$ generators and $\frac{1}{2}(n_c + n_{c_w})^2 + \frac{5}{2}(n_c + n_{c_w})$ constraints, which is a polynomial increase in complexity in terms of both generators and constraints.
\end{remark}

\subsection{Selection of h} \label{sec:selectionofh}

The methods proposed in the previous sections require a choice of $\mathbf{h} \in X = \{\mathbf{G}_x,\mathbf{c}_x,\mathbf{A}_x,\mathbf{b}_x\}$ in order to compute a constrained zonotope enclosure for the prediction step \eqref{eq:prediction0}. As shown in Section \ref{sec:example1}, this choice may drastically affect the accuracy of the obtained enclosure. In the mean value extension for intervals and zonotopes, a usual choice of $\mathbf{h} \in X$ is the center of $X$ \cite{Alamo2005a,Moore2009}. However, with constrained zonotopes, since the center of the CG-rep may not belong to $X$\footnote{From $X = \mbf{c}_x \oplus \mbf{G}_x B_\infty(\mbf{A}_x,\mbf{b}_x)$, $\mathbf{c}_x \notin X$ as long as $\nexists \bm{\xi} \in B_\infty(\mbf{A}_x,\mbf{b}_x)$ satisfying $\mbf{G}_x \bm{\xi} = \bm{0}$.}, a different point $\mathbf{h} \in X$ must be chosen. A simple and inexpensive choice is the center of the interval hull of $X$. Unfortunately, even this point may not belong to $X$ in some cases \footnote{An example is the polytope with vertices $(0,0,0)$, $(1,1,0)$, $(0,1,0)$, $(0,1,1)$.}. Nevertheless, this choice can be applied rigorously by simply checking $\mathbf{h} \in X$ beforehand by solving an LP \cite{Scott2016}. 

In the following, we analyze alternative choices of $\mbf{h} \in X$. Firstly, we focus on suitable choices valid for the mean value extension (Theorem \ref{the:pred2}). This extension relies on the CZ-inclusion (Theorem \ref{the:gzielim}), and therefore requires the computation of a zonotope enclosing $X - \mbf{h}$. In this work, we assume that this zonotope is computed through constraint eliminations (see Remark \ref{rem:contrelim}). Let $\{\mathbf{G}^{(\ell)}, \mathbf{c}^{(\ell)}, \mathbf{A}^{(\ell)}, \mathbf{b}^{(\ell)}\}$ denote the constrained zonotope obtained by reducing to $\ell$ the number of remaining constraints in $X - \mbf{h}$. Following the constraint elimination algorithm in \cite{Scott2016}, for each $\ell=n_c,n_c-1,\dots,1$, the remaining constraints $\mathbf{A}^{(\ell)} \bm{\xi} = \mathbf{b}^{(\ell)}$ are first preconditioned through Gauss-Jordan elimination with full pivoting and then subjected to a rescaling procedure before the next constraint is eliminated. The entire procedure can be represented by the following recursive equations (see Proposition 5 and the Appendix in \cite{Scott2016} for details), where $\bar{(\cdot)}$ denotes variables after preconditioning, $\tilde{(\cdot)}$ denotes variables after rescaling, and $\bm{\Lambda}_\text{G}$, $\bm{\Lambda}_\text{A}$, $\bm{\xi}_\text{m}$, and $\bm{\xi}_\text{r}$ are defined as in \cite{Scott2016}:
\begin{equation} \label{eq:rescaling}
\begin{aligned}
\tilde{\mathbf{c}}^{(\ell)} & = \mathbf{c}^{(\ell)} + \bar{\mathbf{G}}^{(\ell)} \bm{\xi}_\text{m}^{(\ell)}, & \mathbf{c}^{(\ell-1)} & = \tilde{\mathbf{c}}^{(\ell)} + \bm{\Lambda}_\text{G}^{(\ell)} \tilde{\mathbf{b}}^{(\ell)}, \\ \tilde{\mathbf{G}}^{(\ell)} & = \bar{\mathbf{G}}^{(\ell)} \text{diag}(\bm{\xi}_\text{r}^{(\ell)}), & \mathbf{G}^{(\ell-1)} & = \tilde{\mathbf{G}}^{(\ell)} - \bm{\Lambda}_\text{G}^{(\ell)} \tilde{\mathbf{A}}^{(\ell)},\\
\tilde{\mathbf{A}}^{(\ell)} & = \bar{\mathbf{A}}^{(\ell)} \text{diag}(\bm{\xi}_\text{r}^{(\ell)}), & \mathbf{A}^{(\ell-1)} & = \tilde{\mathbf{A}}^{(\ell)} - \bm{\Lambda}_\text{A}^{(\ell)} \tilde{\mathbf{A}}^{(\ell)}, \\
\tilde{\mathbf{b}}^{(\ell)} & = \bar{\mathbf{b}}^{(\ell)} - \bar{\mathbf{A}}^{(\ell)} \bm{\xi}_\text{m}^{(\ell)}, & \mathbf{b}^{(\ell-1)} & = \tilde{\mathbf{b}}^{(\ell)} - \bm{\Lambda}_\text{A}^{(\ell)} \tilde{\mathbf{b}}^{(\ell)}.
\end{aligned}
\end{equation}

Careful examination of the algorithm in \cite{Scott2016} reveals that the actions taken during preconditioning, rescaling, and constraint elimination are all independent of the center of the original constrained zonotope, which in this case is $\mathbf{c}^{(n_c)} = \mathbf{c}_x - \mathbf{h}$. Therefore, with exception of the center, the variables $({\cdot})^{(\ell)}$ can be obtained by eliminating the constraints of $X$ prior to choosing $\mbf{h}$. Considering procedure \eqref{eq:rescaling}, the following corollary provides a choice of $\mbf{h}$ that leads to a tight enclosure by reducing the conservativeness of the CZ-inclusion $\gzinclusion~(\mbf{J},  X - \mbf{h})$.

\begin{corollary} \rm \label{col:C2}
Let $X = \{\mbf{G}_x,\mbf{c}_x,\mbf{A}_x,\mbf{b}_x\} \subset \realset^n$, and consider $\bm{\mu}$, $W$, and $\mbf{J}$ as defined in Theorem \ref{the:pred2}. Assume that $\bar{\mathbf{G}}^{(\ell)}$, $ \bm{\xi}_\text{m}^{(\ell)}$,$\bm{\Lambda}_\text{G}^{(\ell)}$, and $\tilde{\mathbf{b}}^{(\ell)}$ are obtained by eliminating all $n_c$ constraints from $X$ according to \eqref{eq:rescaling} and set
\begin{equation} \label{eq:choice3}
\mathbf{h} = \mathbf{c}_x + \sum_{\ell=1}^{n_c} \left(\bar{\mathbf{G}}^{(\ell)} \bm{\xi}_\text{m}^{(\ell)} + \bm{\Lambda}_\text{G}^{(\ell)} \tilde{\mathbf{b}}^{(\ell)}\right).
\end{equation}
Let $\bar{X} = \{\mbf{G}^{(0)},\mbf{c}^{(0)}\}$ be obtained by eliminating all $n_c$ constraints from $X - \mbf{h}$ according to \eqref{eq:rescaling}, let $\mbf{m} \supseteq (\mbf{J}-\midpoint{\mbf{J}})\mbf{c}^{(0)}$ be computed by standard interval arithmetic, and suppose that $\gzinclusion\left(\mathbf{J},  X - \mathbf{h} \right)$ is computed as in Theorem \ref{the:gzielim} with this choice of $\bar{X}$ and $\mbf{m}$. Finally, let $Z \supseteq \bm{\mu}(\mbf{h},W)$. If $\mbf{h} \in X$, then $\bm{\mu}(X,W) \subseteq  Z \oplus \gzinclusion\left(\mathbf{J},  X - \mathbf{h} \right)$. Moreover, $\gzinclusion\left(\mathbf{J},  X - \mathbf{h} \right) \subseteq \gzinclusion~(\mbf{J},  X - \hat{\mbf{h}})$ for any $\hat{\mbf{h}} \in X$, $\hat{\mbf{h}} \neq \mbf{h}$.
\end{corollary}
\proof 
For $\mbf{h} \in X$, $\bm{\mu}(X,W) \subseteq  Z \oplus \gzinclusion\left(\mathbf{J},  X - \mathbf{h} \right)$ follows directly from Theorem \ref{the:pred2}. Now, let us show that $\gzinclusion\left(\mathbf{J},  X - \mathbf{h} \right) \subseteq \gzinclusion~(\mbf{J},  X - \hat{\mbf{h}})$ holds for any $\hat{\mbf{h}} \in X$, $\hat{\mbf{h}} \neq \mbf{h}$. Recursive computation of \eqref{eq:rescaling} leads to
\begin{equation} \label{eq:pbar}
\mathbf{c}^{(0)} = \mathbf{c}^{(n_c)} + \sum_{\ell=1}^{n_c} \left( \bar{\mathbf{G}}^{(\ell)} \bm{\xi}_\text{m}^{(\ell)} + \bm{\Lambda}_\text{G}^{(\ell)} \tilde{\mathbf{b}}^{(\ell)} \right),
\end{equation}
where $\mathbf{c}^{(n_c)} = \mathbf{c}_x - \mathbf{h}$. Therefore, $\mbf{c}^{(0)} = \bm{0}$ iff $\mbf{h}$ is given by \eqref{eq:choice3}, thus $\mbf{m} = \bm{0}$, and $\text{diam}(\mbf{m}) = \mbf{0}$. Note that in \eqref{eq:czinclusionP}, $\bar{\mbf{M}} \triangleq \mbf{G}^{(0)}$ is invariant with respect to $\mbf{h}$, and since $\mbf{J} \supseteq \nabla_x^T\bm{\mu}(X,W)$, then $\mbf{J}$ is also invariant with respect to $\mbf{h}$. Consequently, the second term in \eqref{eq:czinclusionP} is not a function of $\mbf{h}$. Therefore, $\mbf{P} B_\infty^n \subseteq \hat{\mbf{P}}B_\infty^n = (1/2)\text{diag}(\diam{\hat{\mbf{m}}})B_\infty^n \oplus \mbf{P} B_\infty^n$, with $\hat{\mbf{m}}$ computed using $\hat{\mbf{h}} \in X$. The result then follows from \eqref{eq:czinclusion}.
\qed

By Corollary \ref{col:C2}, the enclosure obtained in Theorem \ref{the:pred2} is tightened by choosing $\mbf{h}$ such that $\mbf{c}^{(0)}$ is equal to zero. Unfortunately, the $\mbf{h}$ given by \eqref{eq:choice3} may not belong to $X$, so an alternative to obtain tight bounds is to reduce the size of the box $\mbf{m}$ by solving 
\begin{equation} \label{eq:choice3optimal}
\underset{\mathbf{h}}{\min}~\{\|\diam{\mathbf{m}}\|_1 : \mathbf{h} \in X\},
\end{equation}
with $\mbf{m} \supseteq (\mbf{J}-\midpoint{\mbf{J}})\bar{\mbf{p}}$ computed using interval arithmetic, where $\bar{\mbf{p}} \triangleq \mbf{c}^{(0)}$. Recall that $\mbf{c}^{(0)}$ is the center of the zonotope obtained by eliminating all the constraints of $X-\mbf{h}$.

\begin{lemma} \rm \label{lem:choosemindiam}
Let $X = \{\mbf{G}_x,\mbf{c}_x,\mbf{A}_x,\mbf{b}_x\} \subset \realset^n$, $\mbf{J} \in \intvalsetmat{n}{n}$. Assume that $\bar{\mathbf{G}}^{(\ell)}$,  $\bm{\xi}_\text{m}^{(\ell)}$, $\bm{\Lambda}_\text{G}^{(\ell)}$, and $\tilde{\mathbf{b}}^{(\ell)}$ are obtained by eliminating all $n_c$ constraints of $X$ according to \eqref{eq:rescaling}. Then $\mbf{h} = \mbf{c}_x + \mbf{G}_x \bm{\xi}^*$ is the solution to \eqref{eq:choice3optimal} iff $\bm{\xi}^*$ is the solution to the linear program
\begin{equation} \label{eq:choosehmindiam}
\underset{\bm{\xi}}{\min}~ \|\bm{\Theta} \bar{\mbf{p}}\|_1, \quad
\text{s.t.} \quad \mathbf{A}_x \bm{\xi} = \mathbf{b}_x, \quad \ninf{\bm{\xi}} \leq 1,
\end{equation}
with $\bar{\mbf{p}} = -\mathbf{G}_x \bm{\xi} + \sum_{\ell=1}^{n_c} \left( \bar{\mathbf{G}}^{(\ell)} \bm{\xi}_\text{m}^{(\ell)} + \bm{\Lambda}_\text{G}^{(\ell)} \tilde{\mathbf{b}}^{(\ell)} \right)$, $\Theta_{jj} = \sum_{i=1}^n \diam{J_{ij}}$, and $\Theta_{ij} = 0$ for $i\neq j$.
\end{lemma}
\proof Each element of $(\mathbf{J} - \midpoint{\mathbf{J}})\in \intvalsetmat{n}{n}$ is a symmetric interval satisfying $(J_{ij} - \midpoint{J_{ij}}) = (1/2) \diam{J_{ij}}[-1,1]$, and for every $a \in \realset$, $a [-1,1] = |a| [-1,1]$ holds. Therefore $m_i = \sum_{j=1}^{n} (1/2) \diam{J_{ij}} |\bar{p}_j| [-1,1]$. Consequently, $\diam{m_i} = \sum_{j=1}^n \diam{J_{ij}} |\bar{p}_j|$, and
\begin{align*}
\|\diam{\mathbf{m}}\|_1 & = \sum_{i=1}^n \sum_{j=1}^n \diam{J_{ij}} |\bar{p}_j| = \sum_{j=1}^n \left( \sum_{i=1}^n \diam{J_{ij}} \right) |\bar{p}_j| \\
 & = \sum_{j=1}^n \Theta_{jj} |\bar{p}_j| =  \|\bm{\Theta} \bar{\mathbf{p}}\|_1.
\end{align*}
The equality $\bar{\mbf{p}} = -\mathbf{G}_x \bm{\xi} + \sum_{\ell=1}^{n_c} \left( \bar{\mathbf{G}}^{(\ell)} \bm{\xi}_\text{m}^{(\ell)} + \bm{\Lambda}_\text{G}^{(\ell)} \tilde{\mathbf{b}}^{(\ell)} \right)$ and the constraints in \eqref{eq:choosehmindiam} follow directly from \eqref{eq:pbar} and $\mbf{h} \in X$. \qed
\color{black}

Lemma \ref{lem:choosemindiam} yields an optimal choice of $\mathbf{h} \in X$ that can be used in Theorem \ref{the:pred2} to reduce conservatism in the CZ-inclusion, and requires only the solution of an LP. Note that formulating \eqref{eq:choosehmindiam} requires the knowledge of $\bar{\mathbf{G}}^{(\ell)}$, $\bm{\xi}_\text{m}^{(\ell)}$, $\bm{\Lambda}_\text{G}^{(\ell)}$, and $\tilde{\mathbf{b}}^{(\ell)}$, which are obtained from the iterated constraint elimination process. As stated before, constraint elimination can be performed over $X$ to obtain the required data prior to the solution of \eqref{eq:choosehmindiam}. Once the optimal $\mathbf{h}$ is obtained, constraint elimination can be repeated, or equivalently, the zonotope obtained using $\mathbf{h} = \bm{0}$ can simply be translated by $-\mathbf{h}$.

\begin{remark} \rm
Note that if the $\mbf{h}$ given by Corollary \ref{col:C2} belongs to $X$, then this coincides with the solution provided by Lemma \ref{lem:choosemindiam}.
\end{remark}

We summarize the proposed choices of $\mathbf{h} \in X$ for use in Theorem \ref{the:pred2} as follows:
\begin{itemize}
	\item[\textbf{\textit{C1)}}] $\mathbf{h}$ is given by the center of the interval hull of $X$, if it satisfies $\mathbf{h} \in X$;
	\item[\textbf{\textit{C2)}}] $\mathbf{h}$ is obtained by solving \eqref{eq:choosehmindiam}.
\end{itemize}

\color{black}

We now focus on suitable choices valid for the first-order Taylor extension. As with the mean value extension, the usual choice of $\mathbf{h} \in X$ in first-order Taylor extensions for intervals and zonotopes is the center of $X$ \cite{Moore2009,Combastel2005}. The next corollary shows that this choice leads to a tight enclosure if it belongs to $X$.

\begin{corollary} \rm \label{col:C3}
Let $Z = \{\mbf{G},\mbf{c},\mbf{A},\mbf{b}\} = X \times W \subset \realset^m$, and consider $\bm{\eta}$, $\tilde{\mbf{c}}$, $\tilde{\mbf{G}}$, $\tilde{\mbf{G}}_\mbf{d}$, $\tilde{\mbf{A}}$, $\tilde{\mbf{b}}$, $\mbf{L}$, and $\mbf{Q}^{[q]}$ as defined in Theorem \ref{thm:firstorder}, with $q \in \{1,2,\ldots,n\}$. If $\mbf{h} = \mbf{c} \in Z$, then $\bm{\eta}(Z) \subseteq \bm{\eta}(\mathbf{h}) \oplus \nabla^T \bm{\eta}(\mathbf{h})(Z - \mathbf{h}) \oplus \tilde{\mathbf{c}} \oplus \left[ \tilde{\mathbf{G}} \,\; \tilde{\mathbf{G}}_{\mathbf{d}} \right] B_\infty(\tilde{\mathbf{A}}, \tilde{\mathbf{b}})$.
\end{corollary}
\proof 
For $\mbf{h} = \mbf{c}$, $\mathbf{L}_{q,:} \supseteq (\mathbf{c} - \mathbf{h})^T \mathbf{Q}^{[q]} = \bm{0}$, $q = 1,2,\ldots,n$. Therefore $\mbf{L} = \mbf{0}$ holds, and $\gzinclusion (\mathbf{L}, (\mathbf{c} - \mathbf{h}) \oplus 2\mathbf{G} B_\infty(\mathbf{A},\mathbf{b}) ) = \mbf{0}$. The result then follows from \eqref{eq:firstorderextension}. 
\qed

By inspecting Corollary \ref{col:C3}, it is clear that the enclosure in \eqref{eq:firstorderextension} is tightened since $\gzinclusion (\mathbf{L}, (\mathbf{c} - \mathbf{h}) \oplus 2\mathbf{G} B_\infty(\mathbf{A},\mathbf{b}) ) = \mbf{0}$. However, since this point may not belong to $X$, a good alternative may be to consider the closest point in $X$ to its center, obtained by means of Proposition \ref{propo:closest}. By the definition of $\mbf{L}$, this heuristic leads to smaller values of $\text{diam}(\mbf{L})$, and therefore reduces the size of $\gzinclusion~(\mathbf{L}, (\mathbf{c} - \mathbf{h}) \oplus 2\mathbf{G} B_\infty(\mathbf{A},\mathbf{b}) )$ (see \eqref{eq:czinclusionP}). A third option is to apply Proposition \ref{propo:center} to obtain an alternative CG-rep of $X$ with any desired center. In this case, the new center is chosen as some point in $X$, $\bar{\mathbf{h}} \in X$, and then $\mathbf{h}$ is chosen as $\mathbf{h} = \bar{\mathbf{h}}$. A simple choice of new center $\bar{\mathbf{h}} \in X$ for such a procedure is the center of the interval hull of $X$. The proposed alternatives for use in Theorem \ref{thm:firstorder} are summarized as follows:
\begin{itemize}
    \item[\textbf{\textit{C3)}}] $\mathbf{h}$ is given by the closest point in $X$ to the center of $X$, computed through Proposition \ref{propo:closest};	
	\item[\textbf{\textit{C4)}}] $\mbf{h}$ is the center of $X$, if it satisfies $\mbf{h} \in X$. Otherwise, $\bar{\mathbf{h}} \in X$ is chosen as the center of the interval hull of $X$, the center of $X$ is moved to $\bar{\mathbf{h}}$ using Proposition \ref{propo:center}, and $\mathbf{h}$ is given by $\mathbf{h}=\bar{\mathbf{h}}$.\footnote{Note that this choice may lead to the same value of $\mbf{h}$ provided by \emph{C1}, but $X$ is described by a different CG-rep with center $\mbf{h}$. If this point is not in $X$, $\bar{\mathbf{h}}$ can be chosen as the point obtained from Proposition \ref{propo:closest} instead.}
\end{itemize}

\subsection{Update step} \label{sec:updatestep}

An enclosure for the prediction step \eqref{eq:prediction0} can be obtained in CG-rep using either Theorem \ref{the:pred2} or Theorem \ref{thm:firstorder}. Therefore, due to linearity of the measurement in \eqref{eq:system}, an exact bound for the update step \eqref{eq:update0} can be directly obtained by computing the generalized intersection of two constrained zonotopes as follows. Given the prediction set $\bar{X}_k$, a constrained zonotope $V$ describing bounds on measurement errors, the current input $\mathbf{u}_k$ and measurement $\mathbf{y}_k$, an exact enclosure for the update step is obtained using the definition \eqref{eq:intersection}, given by 
\begin{equation} \label{eq:update}
\hat{X}_k = \bar{X}_k \cap_{\mathbf{C}} ((\mathbf{y}_k - \mathbf{D}_u \mathbf{u}_k) \oplus (-\mathbf{D}_v V)).
\end{equation}
It is well known that the intersection in \eqref{eq:update} can not be computed exactly using zonotopes, and must be over-approximated \cite{Le2013,Alamo2005a}. As a consequence, the enclosures of the system states obtained after many iterations of prediction and update may be quite conservative using zonotopes. However, with constrained zonotopes all operations in \eqref{eq:update} are easily computed through \eqref{eq:czlimage}--\eqref{eq:czintersection}. These lead to an enclosure with $n_g + n_{g_v}$ generators, and $n_c + n_{c_v} + n_y$ constraints, where $n_g$ and $n_c$ are the number of generators and constraints of $\bar{X}_k$, respectively

\begin{remark} \rm

Iterated computations of the proposed extensions (Theorems \ref{the:pred2} and \ref{thm:firstorder}) and \eqref{eq:update} result in at most a quadratic increase in the complexity of the CG-rep \eqref{eq:czdefinition}. As with zonotopes, this can be effectively addressed using \emph{order reduction} algorithms \cite{Kopetzki2017,Scott2018} that over-approximate a constrained zonotope by another with lower complexity. Efficient methods for reducing the number of generators and constraints of the CG-rep \eqref{eq:czdefinition} with reasonable conservativeness were proposed in \cite{Scott2016}.

\end{remark}

\subsection{Complexity analysis}

Table \ref{tab:complexityestimators} shows the computational complexity\footnote{We use the standard $O({\cdot})$ notation defined in \citep{Cormen2009}.} of our methods for the prediction and update steps, as well as complexity reduction to the same number of generators and constraints of the set prior to prediction. Specifically, we use the mean value extension (Theorem \ref{the:pred2}) and the first-order Taylor extension (Theorem \ref{thm:firstorder}) for the prediction steps, while the update steps are both given by the generalized intersection \eqref{eq:update}. These methods are denoted by CZMV and CZFO, respectively. The computational complexities of their zonotope counterparts are also presented for comparison, denoted analogously by ZMV and ZFO, which use the mean value approach in \cite{Alamo2005a} and the first-order Taylor approach in \cite{Combastel2005}, respectively, for the prediction step. The update algorithm proposed in \cite{Bravo2006} is used for both ZMV and ZFO because it provided the best trade-off between accuracy and efficiency in our numerical experiments with zonotopes. Complexity reduction is applied after the update step in all four methods using the reduction methods in \cite{Scott2016} for constrained zonotopes and Method 4 in \cite{Scott2018} for zonotopes. For constrained zonotopes, constraint elimination is performed prior to generator reduction. The complexities in Table \ref{tab:complexityestimators} take into account the growth of the number of generators and constraints after each step (see Remarks \ref{rem:meanvaluecomplexity} and \ref{rem:firstordercomplexity}). The dimensions in Table \ref{tab:complexityestimators} are specified by the definitions $\hat{X}_{k-1} = \{\mbf{G}_x,\mbf{c}_x\}$, $W = \{\mbf{G}_w, \mbf{c}_w\}$, and $V = \{\mbf{G}_v, \mbf{c}_v\}$ or $\hat{X}_{k-1} = \{\mbf{G}_x,\mbf{c}_x,\mbf{A}_x,\mbf{b}_x\}$, $W = \{\mbf{G}_w, \mbf{c}_w, \mbf{A}_w, \mbf{b}_w\}$, and $V = \{\mbf{G}_v, \mbf{c}_v, \mbf{A}_v, \mbf{b}_v\}$ with $\mbf{G}_x \in \realsetmat{n}{n_g}$, $\mbf{c}_x \in \realset^n$, $\mbf{A}_x \in \realsetmat{n_c}{n_g}$, $\mbf{b}_x \in \realset^{n_c}$, $\mbf{G}_w \in \realsetmat{n_w}{n_{g_w}}$, $\mbf{c}_w \in \realset^{n_w}$, $\mbf{A}_w \in \realsetmat{n_{c_w}}{n_{g_w}}$, $\mbf{b}_w \in \realset^{n_{c_w}}$, $\mbf{G}_v \in \realsetmat{n_v}{n_{g_v}}$, $\mbf{c}_v \in \realset^{n_v}$, $\mbf{A}_v \in \realsetmat{n_{c_v}}{n_{g_v}}$, $\mbf{b}_v \in \realset^{n_{c_v}}$, $\mbf{u}_k \in \realset^{n_u}$, and $\mbf{y}_k \in \realset^{n_y}$. For simplicity, we define $m = n + n_w$, $m_g = n_g + n_{g_w}$, $m_c = n_c + n_{c_w}$, $\delta_n = n - n_y$, $\delta_w = n_{g_w} - n_{c_w}$, $\delta_v = n_{g_v} - n_{c_v}$, and $\tilde{\delta} = m_g^2 - m_c^2$. Moreover, we consider that scalar real function and scalar inclusion function evaluations have complexity $O(1)$. These correspond to evaluations of the nonlinear dynamics in \eqref{eq:system} and its derivatives using real and interval arithmetic, respectively. The complexities of the basic operations on zonotopes and constrained zonotopes used to derive the figures in Table \ref{tab:complexityestimators} can be found in the Appendix. 

The dominant terms in the prediction step of CZMV come from the computation of the interval hulls of $X$ and $W$ and the CZ-inclusions $\gzinclusion\left(\mathbf{J},  X - \mathbf{h} \right)$ and $\gzinclusion\left(\mathbf{J}_w,  W - \mathbf{h}_w \right)$ in Theorem \ref{the:pred2} and Remark \ref{rem:affine}. In the case of CZFO, the dominant terms come from the computation of the interval matrices $\tilde{\mbf{Q}}^{[q]}$, the interval hull of $Z = X \times W$, and the CZ-inclusion $\gzinclusion (\mathbf{L}, (\mathbf{c} - \mathbf{h}) \oplus 2\mathbf{G} B_\infty(\mathbf{A},\mathbf{b}) )$ in Theorem \ref{thm:firstorder}. Note that the worst-case complexities of the prediction steps of our methods are higher than the zonotope methods, while the update steps are cheaper due to the generalized intersection \eqref{eq:update}. Even so, the complexity of the proposed methods are still polynomial. For a simplified analysis, assuming that all of the variables in Table \ref{tab:complexityestimators} increase linearly with $n$, the total complexities for ZMV, ZFO, CZMV and CZFO are $O(n^4)$, $O(n^5)$, $O(n^5)$, and $O(n^8)$, respectively. On the other hand, even basic polytope operations are known to be exponential \cite{Hagemann2015}. Besides, despite the higher complexities of CZMV and CZFO in comparison to the zonotope methods, they provide more accurate enclosures as shown in the next section. 

\begin{table*}[!tb]
  \scriptsize
    \centering
	\caption{Computational complexity $O(\cdot)$ of the state estimators.}
	\begin{tabular}{c c c} \hline
		Step & ZMV & ZFO \\ \hline
	    Prediction & $n^2n_g + nn_wn_{g_w}$ & $n(m^2m_g + mm_g^2)$ \\
	    Update & $n_y(n^3(m_g+n)+n^2(m_g+n)^2+n_u+n_vn_{g_v})$ & $n_y(n^3(m_g^2+n)+n^2(m_g^2+n)^2+n_u+n_vn_{g_v})$  \\
		Reduction & $n^2(m_g+n) + n(n_{g_w}+n)(m_g+n)$ & $n^2(m_g^2+n) + n(m_g^2+n)^2$ \\
		\hline
		Step & CZMV & CZFO \\ \hline
		\multirow{1}{*}{Prediction} & $n^2n_g + nn_wn_{g_w} + (nn_g+n_c)(n_g+n_c)^3 + (n_wn_{g_w}+n_{c_w})(n_{g_w}+n_{c_w})^3$ & $n(m^2m_g+mm_g^2) + (mm_g+m_c)(m_g+m_c)^3$\\
		Update & $n_yn(m_g+n) + n_yn_u + n_yn_vn_{g_v}$ & $n_yn(m_g^2+n) + n_yn_u + n_yn_vn_{g_v}$ \\
		\multirow{2}{*}{Reduction} & $(n_{c_w}+n_{c_v}+n_y)(m_g+m_c+n_{g_v}+n_{c_v}+n+n_y)^3$ & $(m_c^2+n_{c_v}+n_y)(m_g^2+m_c^2+n_{g_v}+n_{c_v}+n+n_y)^3 $ \\
		& $+ (n+n_c)^2(n_g+\delta_n+\delta_w+\delta_v)+(n+n_c)(\delta_n+\delta_w+\delta_v)(n_g+\delta_n+\delta_w+\delta_v)$ & $+ (n+n_c)^2(\tilde{\delta}+\delta_n + \delta_v)+(n+n_c)(\tilde{\delta}+\delta_n+\delta_v)^2$ \\
		\hline		
	\end{tabular} \normalsize
 	\label{tab:complexityestimators}
\end{table*}

\section{Numerical examples} \label{sec:numericalexamples}


This section presents numerical results for the two new set-valued state estimation methods enabled by the results in the previous section. The imposed limits on the complexity of the sets used are described separately for each example below.

\subsection{Example 1} \label{sec:example1}

To demonstrate the effect of the different choices of $\mathbf{h}$, we first analyze one iteration of the prediction step for the nonlinear system \cite{Raimondo2012}
\begin{equation} \label{eq:model1}
\begin{aligned}
x_{1,k} & = 3 x_{1,k-1} - \frac{x_{1,k-1}^2}{7} - \frac{4 x_{1,k-1} x_{2,k-1}}{4 + x_{1,k-1}} + w_{1,k-1}, \\
x_{2,k} & = -2 x_{2,k-1} + \frac{3 x_{1,k-1} x_{2,k-1}}{4 + x_{1,k-1}} + w_{2,k-1},
\end{aligned}
\end{equation}
with
\begin{equation} \label{eq:reachX0}
	 X_0 = \left\{ \begin{bmatrix} 0.2 & 0.4 & 0.2 \\ 0.2 & 0 & -0.2 \end{bmatrix}, \begin{bmatrix} -1 \\ 1 \end{bmatrix}, \begin{bmatrix} 2 & 2 & 2  \end{bmatrix}, -3 \right\},
\end{equation}
where $\mathbf{w}_k \in \realset^2$ denotes process uncertainties, which are zero in this first scenario.

\begin{figure}[!tb]
		\centering{
			\def\svgwidth{\columnwidth}
			{\scriptsize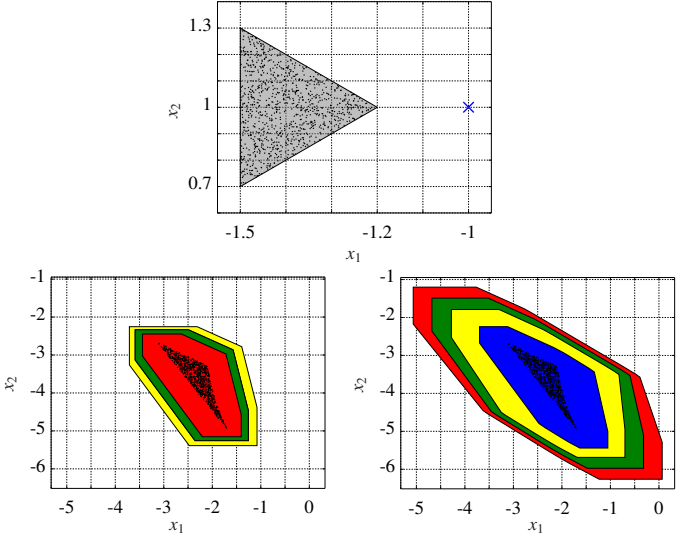}
            \caption{Top: The constrained zonotope $X_0$ with `$\times$' denoting its center. Left: Enclosures obtained by applying Theorem \ref{the:pred2} to \eqref{eq:model1}. Right: Enclosures obtained by applying Theorem \ref{thm:firstorder} to \eqref{eq:model1}. The real vector $\mathbf{h}$ is determined by \emph{C1} (green), \emph{C2} (red), \emph{C3} (yellow), and \emph{C4} (blue). For the mean value extension (left), \emph{C4} is overlapped with \emph{C1}. Black dots denote uniform samples from $X_0$ propagated through \eqref{eq:model1}.}\label{fig:reachable}}
\end{figure}

Figure \ref{fig:reachable} shows the constrained zonotope $X_0$ and the enclosures of the one-step reachable set obtained by Theorem \ref{the:pred2} using \emph{C1}--\emph{C4}. Since the complexity of the enclosure for \emph{C4} is higher than for the other methods (see Proposition \ref{propo:center}), the reduction methods in \cite{Scott2016} were used to reduce the number of generators and constraints in this enclosure to match the other methods before comparison. In this example, the choice of $\mathbf{h}$ has a moderate impact in the enclosure computed by Theorem \ref{the:pred2}, with \emph{C2} providing the least conservative result, as expected. Therefore, \emph{C2} is employed in Theorem \ref{the:pred2} henceforth. 

Figure \ref{fig:reachable} also shows the enclosures of the one-step reachable set obtained by Theorem \ref{thm:firstorder} with \emph{C1}--\emph{C4}. Clearly, the enclosures produced by Theorem \ref{thm:firstorder} are strongly affected by the choice of $\mathbf{h}$, with \emph{C4} providing the least conservative result. In addition, note that the enclosures provided by the first-order Taylor extension are more conservative than those obtained by the mean value extension. However, experience with zonotopes and intervals (see \cite{Raimondo2012} for detailed examples) suggests that the relative merits of these two methods will depend on the dynamics of the system, as well as the shape and size of the set $X_0$, and the maximum allowed number of generators and constraints. This is corroborated by the next results. 

We consider now the linear measurement equation
\begin{equation} \label{eq:model1meas}
\begin{bmatrix} y_{1,k} \\ y_{2,k} \end{bmatrix} = \begin{bmatrix} 1 & 0 \\ -1 & 1 \end{bmatrix} \begin{bmatrix} x_{1,k} \\ x_{2,k} \end{bmatrix} + \begin{bmatrix} v_{1,k} \\ v_{2,k} \end{bmatrix},
\end{equation}
with bounds $\ninf{\mathbf{w}_k} \leq 0.4$ and $\ninf{\mathbf{v}_k} \leq 0.4$, where $\mathbf{v}_k \in \realset^2$ denote measurement uncertainties. The initial states $\mathbf{x}_0$ are bounded by the zonotope\footnote{Note that $X_0$, $W$ and $V$ are expressed as zonotopes for the purpose of a fair comparison with the zonotope methods.}
\begin{equation} \label{eq:model1X0}
X_0 = \left\{ \begin{bmatrix} 0.1 & 0.2 & -0.1 \\ 0.1 & 0.1 & 0 \end{bmatrix}, \begin{bmatrix} 0.5 \\ 0.5 \end{bmatrix}\right\}.
\end{equation}

\begin{figure}[!tb]
	\centering{
		\def\svgwidth{0.8\columnwidth}
		{\footnotesize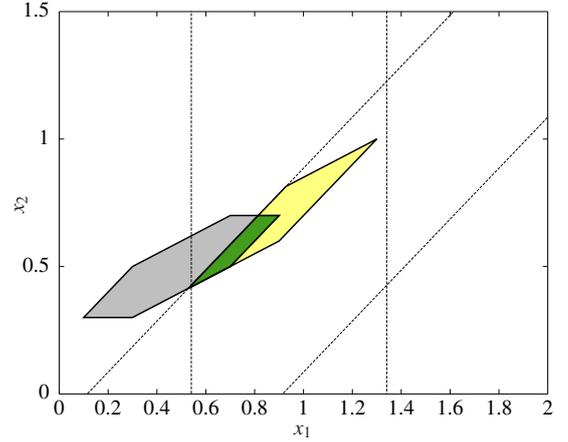}
		\caption{The initial set $X_0$ (gray), the initial bounded uncertain measurements (dashed lines), the intersection computed as in \cite{Bravo2006} (yellow), and the constrained zonotope (green) computed by \eqref{eq:update}.}\label{fig:model1_update}}
\end{figure}

To generate process measurements, \eqref{eq:model1} was simulated with $\mathbf{x}_0 = (0.8,0.65) \in X_0$ and process and measurement uncertainties drawn from uniform random distributions. The number of generators and constraints of the constrained zonotopes was limited to 20 and 5, respectively, while the number of generators of the zonotopes was limited to 20. Figure \ref{fig:model1_update} shows the results of the initial update step using the intersection algorithm in \cite{Bravo2006} and the generalized intersection \eqref{eq:update} computed using \eqref{eq:czintersection}, which yields a constrained zonotope. Clearly, since the generalized intersection is not a symmetric set, it cannot be described by a zonotope. In contrast, the resulting constrained zonotope corresponds to the exact intersection, providing far less conservative bounds in the first update step.

\begin{figure}[!tb]
	\centering{
		\def\svgwidth{\columnwidth}
		{\footnotesize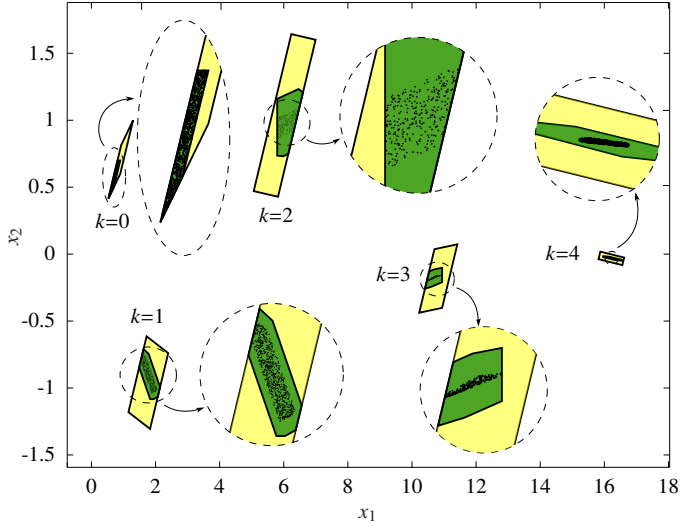}
		\caption{Results from the first four time steps of set-valued state estimation (after update) using the constrained zonotopic method CZMV (green) and the zonotopic method {\alamobravo} (yellow). Black dots denote uniform samples from $X_0$ propagated through \eqref{eq:model1} that are consistent with the current measurement.} \label{fig:model1_meanvalue_4iter}}
\end{figure}

Figure \ref{fig:model1_meanvalue_4iter} shows the first four time steps of CZMV with $\mathbf{h}$ given by \emph{C2} in a scenario without process uncertainties ($\mathbf{w}_k = \bm{0}$). For comparison, the zonotopes computed using ZMV are also depicted. CZMV provides much less conservative enclosures than {\alamobravo} for this example, demonstrating the effectiveness of the proposed nonlinear state estimation strategy. 

\begin{figure}[!tb]
	\centering{
		\def\svgwidth{\columnwidth}
		{\scriptsize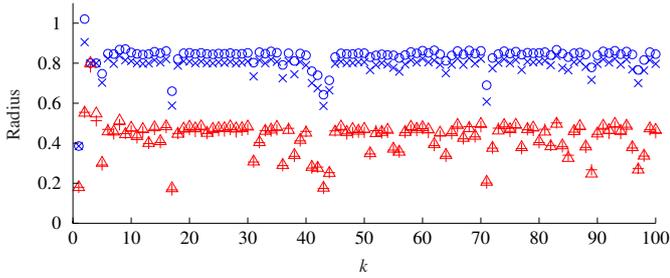}
		\caption{Radii of the update sets obtained by applying CZMV ($\triangle$), {\alamobravo} ($\circ$), CZFO ($+$), and ZFO ($\times$) to \eqref{eq:model1} with process disturbances.}\label{fig:model1_comparison_radiusw}}
\end{figure}

Figure \ref{fig:model1_comparison_radiusw} shows the radii (half the length of the longest edge of the interval hull) of the sets provided by CZMV using \emph{C2} and {\alamobravo} over 100 time steps considering process disturbances. Since \eqref{eq:model1} is affine in $\mathbf{w}_k$, the enclosure $Z \supseteq \bm{\mu}(\mathbf{h},W)$ in Theorem \ref{the:pred2} was computed as described at the end of Remark \ref{rem:affine}. CZMV provided less conservative bounds than the zonotopes computed by {\alamobravo}, with a CZMV-to-{\alamobravo} average radius ratio (ARR, i.e., the ratio of the radius of the CZMV set at $k$ over the ZMV set at $k$ averaged over all time steps $k$) of only $51.4\%$.  
%
Figure \ref{fig:model1_comparison_radiusw} also compares the radii of the update sets computed by {\combastelbravo} and CZFO with $\mathbf{h}$ given by \emph{C4}. As in the previous case, CZFO provides less conservative bounds than {\combastelbravo}, with the CZFO-to-ZFO ARR being only $53.66\%$. The size of the sets provided by CZMV and CZFO were quite similar, with CZFO being less conservative (the CZFO-to-CZMV ARR was $98.75\%$). In both experiments, the number of generators and constraints were limited to 20 and 5, respectively. The ARR for different numbers of constraints are shown in Table \ref{tab:example2ARR}, with the average computed considering in addition simulations with different numbers of generators. Execution times are shown in Table \ref{tab:example2times}. These were obtained using MATLAB 9.1 with CPLEX 12.8 and INTLAB 9, in a laptop with 8GB RAM and an Intel Core i7 4510U 3.1 GHz processor.

The use of constrained zonotopes in CZMV and CZFO results in sets that are slightly more complex than those generated by ZMV and ZFO (specifically, the set description involves five equality constraints that are not present in the zonotopes from ZMV and ZFO). However, this example shows that this increase in complexity is compensated by greatly improved accuracy.

\begin{table}[!htb]
	\footnotesize
	\centering
	\caption{
	Average radius ratio of the estimators with varying numbers of constraints. Each average is taken over 3 separate simulations using $n_g \in \{ 8,12,20\}$.}
	\begin{tabular}{c c c c} \hline
		$n_c$ & CZMV/ZMV & CZFO/ZFO & CZFO/CZMV \\ \hline
	    $1$ & $54.1\%$ & $59.8\%$ & $104.5\%$ \\	    
	    $3$ & $51.6\%$ & $54.0\%$ & $99.0\%$ \\		
		$5$ & $51.6\%$ & $53.7\%$ & $98.5\%$ \\
		\hline
	\end{tabular} \normalsize
	\label{tab:example2ARR}
\end{table}

\begin{table}[!htb]
	\footnotesize
	\centering
	\caption{Average total times per iteration of the estimators with varying numbers of constraints. Each average is taken over 30 separate simulations using $n_g \in \{8,12,20\}$. Times spent only on complexity reduction are shown in parenthesis.}
	\begin{tabular}{c c c c c} \hline
		$n_c$ & ZMV & ZFO & CZMV & CZFO \\ \hline
		$0$ & $5.79~(0.24)$ ms & $14.03~(2.93)$ ms & -- & -- \\
		$1$ & -- & -- & $19.2~(1.6)$ ms & $76.4~(57.1)$ ms \\
	    $3$ & -- & -- & $20.7~(1.8)$ ms & $4.28~(4.26)$ s \\		
	    $5$ & -- & -- & $22.5~(1.9)$ ms & $8.93~(8.91)$ s \\
		\hline
	\end{tabular} \normalsize
	\label{tab:example2times}
\end{table}

\subsection{Example 2}
Consider the quadrotor unmanned aerial vehicle (UAV) described in \cite{Mistler2001} with state vector $\bm{\zeta} = [x\,\;y\,\;z\,\;u\,\;v\,\;w\,\;\phi\,\;\theta\,\;\psi\,\;p\,\;q\,\;r]^T$, where $[x\,\;y\,\;z]^T$ is the position of the UAV with respect to the inertial frame, $[u\,\;v\,\;w]^T$ is the velocity vector expressed in the inertial frame, $[\phi\,\;\theta\,\;\psi]^T$ are Euler angles describing the orientation of the UAV, and $[p\,\;q\,\;r]^T$ is the angular velocity vector expressed in the body frame. The equations of motion are \cite{Mistler2001}:
\begin{equation} \label{eq:quadrotor}
\dot{\bm{\zeta}} = 
\left\{ \begin{aligned}
\dot{x} & = u, \\
\dot{y} & = v, \\
\dot{z} & = w, \\
\dot{u} & = \frac{1}{m} (\cos\psi \sin\theta \cos\phi + \sin\psi \sin\phi)U_1 + \frac{1}{m}D_x, \\
\dot{v} & = \frac{1}{m} (\sin\psi \sin\theta \cos\phi - \cos\psi \sin\phi)U_1 + \frac{1}{m}D_y, \\
\dot{w} & = -g + \frac{1}{m} (\cos\theta \cos\phi)U_1 + \frac{1}{m}D_z, \\
\dot{\phi} & = p + q \sin\phi \tan\theta + r \cos\phi \tan\theta, \\
\dot{\theta} & = q \cos\phi- r \sin\phi, \\
\dot{\psi} & = q \sin\phi \sec\theta + r \cos\phi \sec\theta, \\
\dot{p} & = \frac{I_{yy} - I_{zz}}{I_{xx}} qr + \frac{l}{I_{xx}} U_2, \\
\dot{q} & = \frac{I_{zz} - I_{xx}}{I_{yy}} pr + \frac{l}{I_{yy}} U_3, \\
\dot{r} & = \frac{I_{xx} - I_{yy}}{I_{zz}} pq + \frac{1}{I_{zz}} U_4,
\end{aligned} \right.
\end{equation}
where $m$, $I_{xx}$, $I_{yy}$, $I_{zz}$, and $l$ are physical parameters, $g$ is the gravitational acceleration, $U_1$ is the total thrust generated by the propellers, $U_2$ is the difference of thrusts between the left and right propellers, $U_3$ is the difference of thrusts between the front and back propellers, $U_4$ is the difference of torques between clockwise and counter-clockwise turning propellers, and $\mathbf{d} = [D_x \,\; D_y \,\; D_z]^T$ are disturbance forces applied to the UAV with $\ninf{\mathbf{d}} \leq 1$. The experiment consists in obtaining guaranteed bounds on the system states $\bm{\zeta}$ while the quadrotor UAV tracks a vertical helix trajectory defined by the reference values $x^\text{ref}(t) = \frac{1}{2} \cos\left(\frac{t}{2}\right)$, $y^\text{ref}(t) = \frac{1}{2} \sin\left(\frac{t}{2}\right)$, $z^\text{ref}(t) = 1 + \frac{t}{10}$, $\psi^\text{ref} = \frac{\pi}{3}$,
subject to the disturbance forces described by $D_x = 1$ N for $t \in [5,15)$ s, $D_y = 1$ N for $t \in [8,15)$ s, and $D_z = 1$ N for $t \in [10,15)$ s. These forces are zero otherwise.

The dynamic feedback controller in \cite{Mistler2001} is used to track the reference trajectory above\footnote{In this experiment, the control action is computed using the real states $\bm{\zeta}_k$. The approach in \cite{Rego2018b} can be used for feedback connection using a point that belongs to $\hat{X}_k$.}. The simulation parameters are $m = 0.7$ Kg, $l = 0.3$ m, $I_{xx} = I_{yy} = I_{zz} = 1.2416$ Kg${\cdot}$m$^2$, $g = 9.81$ m/s$^2$. We consider a realistic scenario in which the available measurements are provided by sensors located at the quadrotor UAV, which include: (i) a Global Positioning System (GPS); (ii) a barometer; and (iii) an Inertial Measurement Unit (IMU). The measurements are affected by bounded uncertainties as described in Table \ref{tab:results_sensorparameters}. The velocity vector $[u\,\;v\,\;w]^T$ is not measured.


\begin{table}[!htb]
	\footnotesize
	\centering
	\caption{Measured variables with error bounds.}
	\begin{tabular}{c c c} \hline
		Sensor & Variables & Noise bounds \\ \hline
		GPS & $\{x,y\}$ & $\pm\!\!$ $0.15\!$ m \\ 
		Barometer & $\{z\}$ & $\pm\!\!$ $0.51\!$ m \\ 
		\multirow{2}{*}{IMU} & $\{\phi,\theta,\psi\}$ & $\pm\!\!$ $2.618 {\cdot} 10^{-3}\!$ rad \\
		& $\{p,q,r\}$ & $\pm\!\!$ $16.558 {\cdot} 10^{-3}\!$ rad/s \\
		\hline
	\end{tabular} \normalsize
	\label{tab:results_sensorparameters}
\end{table}

\begin{figure}[!tb]
	\centering{
		\def\svgwidth{\columnwidth}
		{\scriptsize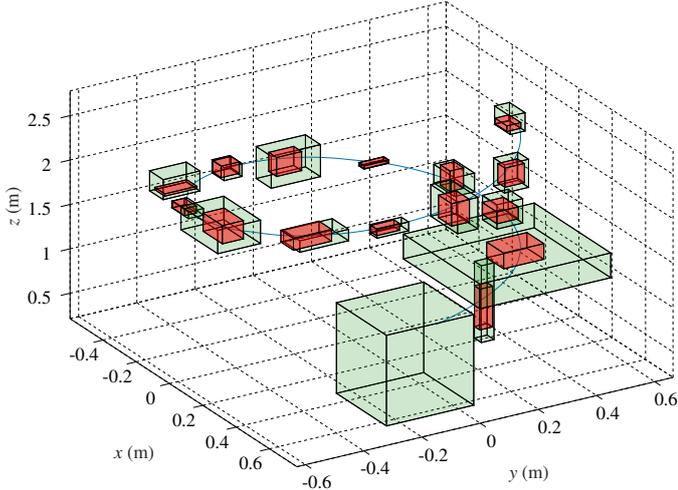}
		\caption{The trajectory performed by the quadrotor UAV (solid line) and the interval hulls of the constrained zonotopes (red boxes) and zonotopes (green boxes) estimated by CZMV and ZMV, respectively, projected onto $(x, y, z)$.}\label{fig:quadrotor_trajectory_boxes}}
\end{figure}

The nonlinear equations \eqref{eq:quadrotor} were discretized by Euler approximation with sampling time $0.01$ s. The initial states $\bm{\zeta}_0$ are bounded by $X_0 = \{\mathbf{G}_0,\bm{0}\}$, where $\mathbf{G}_0 = \text{diag}\left(2, 2, 2, 1, 1, 1, \frac{\pi}{6}, \frac{\pi}{6}, \frac{\pi}{2}, \frac{\pi}{12}, \frac{\pi}{12}, \frac{\pi}{12}\right).$
%
%
To generate process measurements, the discrete-time dynamics were simulated with $\bm{\zeta}_0  = [0.5\,\;0\,\;1\,\; \zeros{1}{5} \,\; \pi/3\,\; \zeros{1}{3} ]^T \in X_0$ and process and measurement noises drawn from uniform distributions. Figure \ref{fig:quadrotor_trajectory_boxes} shows the trajectory performed by the quadrotor UAV along with the interval hulls\footnote{The conversion from CG-rep to H-rep (see \cite{Scott2016}) for the purposes of exact drawing is intractable for the constrained zonotopes in this example.} of the enclosures computed by the methods CZMV and {\alamobravo}, projected onto $(x,y,z)$-axes. CZMV was implemented with $\mathbf{h}$ given by \emph{C2}, and since \eqref{eq:quadrotor} is also affine in $\mathbf{w}_k \triangleq \mathbf{d}_k$, Theorem \ref{thm:firstorder} was implemented with $Z \supseteq \bm{\mu}(\mathbf{h},W)$ computed as described at the end of Remark \ref{rem:affine}. The number of constraints and generators of the computed constrained zonotopes was limited to 40 and 12, respectively, while the number of generators of the computed zonotopes was limited to 40.  

\begin{figure}[!tb]
	\centering{
		\def\svgwidth{\columnwidth}
		{\scriptsize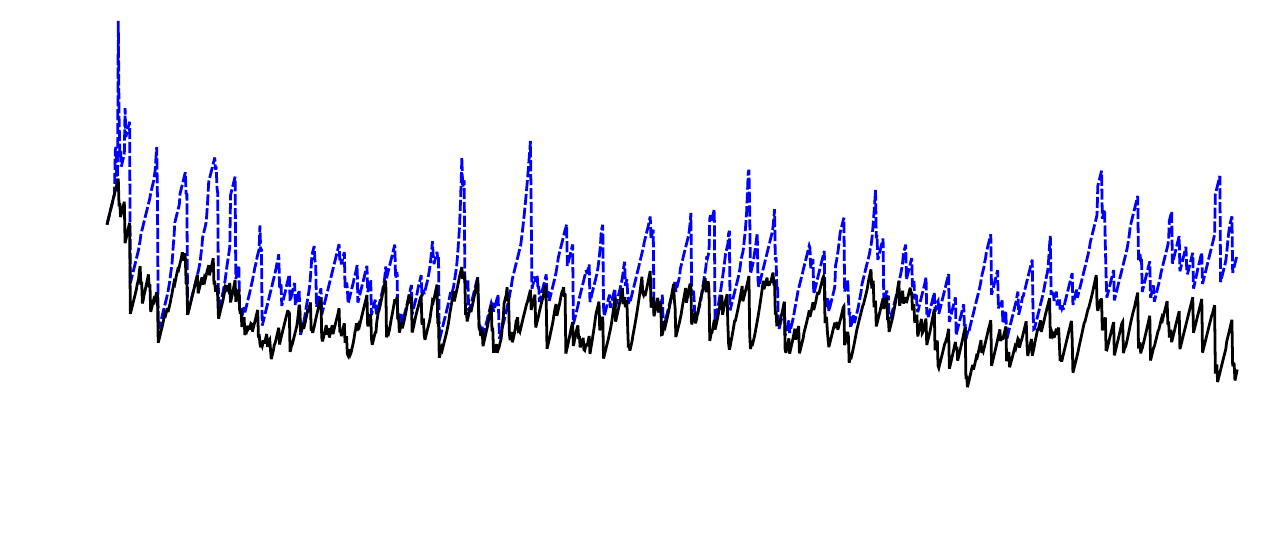}
		\caption{Radii of the update sets computed by CZMV ($-$) and {\alamobravo} ($--$) for the quadrotor UAV experiment.}\label{fig:quadrotor_meanvalue_radius}}
\end{figure}

\begin{figure}[!tb]
	\centering{
		\def\svgwidth{\columnwidth}
		{\scriptsize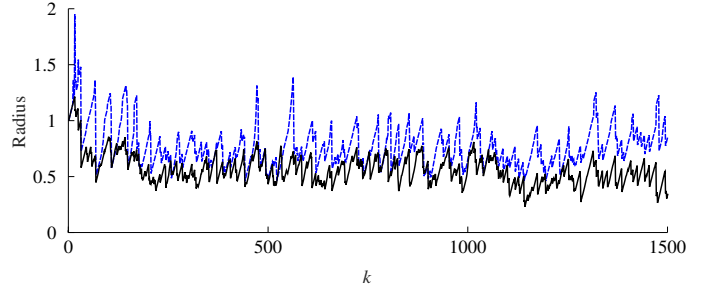}
		\caption{Radii of the update sets computed by CZFO ($-$) and {\combastelbravo} ($--$) for the quadrotor UAV experiment.}\label{fig:quadrotor_firstorder_radius}}
\end{figure}

The interval hulls of the constrained zonotopes obtained by CZMV were smaller than those from {\alamobravo}, demonstrating the accuracy of the proposed method. Figure \ref{fig:quadrotor_meanvalue_radius} shows the radii of the constrained zonotopes and zonotopes computed by CZMV and {\alamobravo}, respectively. Both algorithms were capable of providing tight bounds on the system states $\bm{\zeta}_k \in \realset^{12}$. Nevertheless, CZMV provided less conservative bounds than {\alamobravo}, even for a high-order nonlinear dynamical system such as \eqref{eq:quadrotor} (the CZMV-to-{\alamobravo} ARR was $74.41\%$). %
Finally, Figure \ref{fig:quadrotor_firstorder_radius} compares the radii of the update sets computed by {\combastelbravo} and CZFO with $\mathbf{h}$ given by \emph{C3}\footnote{The increased complexity of the constrained zonotopes provided by \emph{C4} proved to be intractable for this example.}. Once again, CZFO provided less conservative bounds than {\combastelbravo} (the CZFO-to-{\combastelbravo} ARR was $74.45\%$). The results from CZMV and CZFO were again very similar, with CZMV providing marginally better results (the CZMV-to-CZFO ARR was $99.93\%$). The ARR for different numbers of constraints are shown in Table \ref{tab:example3ARR}, with the average computed considering in addition simulations with different numbers of generators. Execution times are shown in Table \ref{tab:example3times}. Note that most of the CZFO execution times are smaller than the ones presented in Table \ref{tab:example2times}. This might be counter-intuitive since the current example has more state variables. However the use of \emph{C4} in Example 1 results in a relatively more complex enclosure and therefore requires a much higher execution time for generator reduction and constraint elimination. Note that in this example, the computational times of the state estimators were greater than the considered sampling time of $0.01$ s. Nevertheless, this fact does not invalidate the obtained results since these were run in a numerical simulation. Better times can be achieved by optimized implementation of the algorithms and using more powerful hardware, for instance. Besides, note that even though the current execution times of CZMV and CZFO would in principle prevent their use in fast applications, the improved accuracy can significantly reduce the number of time steps required for guaranteed fault detection and isolation for systems in which execution time is not critical.

\begin{table}[!htb]
	\footnotesize
	\centering
	\caption{Average radius ratio of the estimators with varying numbers of constraints. Each average is taken over 3 separate simulations with $n_g \in \{ 25,30,40\}$.}
	\begin{tabular}{c c c c} \hline
		$n_c$ & CZMV/ZMV & CZFO/ZFO & CZFO/CZMV \\ \hline
	    $3$ & $82.5\%$ & $82.2\%$ & $99.9\%$ \\
	    $6$ & $76.7\%$ & $77.1\%$ & $100.7\%$ \\		
		$12$ & $75.0\%$ & $74.8\%$ & $99.7\%$ \\
		\hline
	\end{tabular} \normalsize
	\label{tab:example3ARR}
\end{table}

\begin{table}[!htb]
	\footnotesize
	\centering
	\caption{Average total times per iteration of the estimators with varying numbers of constraints. Each average is taken over 15 separate simulations using $n_g \in \{25,30,40\}$. Times spent only on complexity reduction are shown in parenthesis.}
	\begin{tabular}{c c c c c} \hline
		$n_c$ & ZMV & ZFO & CZMV & CZFO \\ \hline
		$0$ & $53.3~(1.0)$ ms & $174.2~(54.2)$ ms & -- & -- \\
		$3$ & -- & -- & $104.7~(7.3)$ ms & $266.9~(128.6)$ ms \\
	    $6$ & -- & -- & $109.5~(8.4)$ ms & $615.0~(471.6)$ ms \\		
		$12$ & -- & -- & $127.8~(12.3)$ ms & $2.62~(2.46)$ s \\
		\hline
	\end{tabular} \normalsize
	\label{tab:example3times}
\end{table}

\section{Conclusions and future work} \label{sec:conclusion}

This paper proposed two novel approaches for set-valued state estimation of nonlinear discrete-time systems with unknown-but-bounded process and measurement uncertainties. A mean value extension and a first-order Taylor extension were developed based on constrained zonotopes, a generalization of zonotopes capable of describing strongly asymmetric convex sets. In addition, measurement data were effectively taken into account by means of generalized intersection, which resulted in far less conservative results than existing methods based on zonotopes. The accuracy of the proposed methods was demonstrated by means of three numerical examples, the third one being an experiment with a quadrotor UAV, considering a realistic scenario with uncertain measurements provided by sensors located on the aircraft. In the latter, execution times were longer than the considered sampling time. Nevertheless, an optimized implementation of the methods, as well as more powerful hardware, implementation in C++ and parallelization techniques, could be used to achieve better times. This issue is left as a future work seeking the practical implementation in a real aircraft.


\section*{Appendix. Computational complexity details}

Table \ref{tab:complexity} shows the computational complexity of the basic operations used in the zonotope and constrained zonotope methods. These complexities assume generic inputs with dimensions $\mbf{R} \in \realsetmat{n_r}{n}$ in \eqref{eq:limage}; $Y =\{\mbf{G}_y, \mbf{c}_y, \mbf{A}_y, \mbf{b}_y\}$ in \eqref{eq:intersection}, with $\mbf{G}_y \in \realsetmat{n_r}{n_{g_r}}$, $\mbf{c}_y \in \realset^{n_r}$, $\mbf{A}_y \in \realsetmat{n_{c_r}}{n_{g_r}}$, and $\mbf{b}_y \in \realset^{n_{c_r}}$; $\mbf{J} \in \intvalsetmat{n_s}{n}$ in Theorem \ref{the:gzielim}; $X = \{\mbf{G}_x,\mbf{c}_x\}$ or $X = \{\mbf{G}_x,\mbf{c}_x,\mbf{A}_x,\mbf{b}_x\}$, with $\mbf{G}_x \in \realsetmat{n}{n_g}$, $\mbf{c}_x \in \realset^n$, $\mbf{A}_x \in \realsetmat{n_c}{n_g}$, and $\mbf{b}_x \in \realset^{n_c}$; $k_g$ and $k_c$ are the number of generators and constraints removed in the order reduction process, respectively. `Set inclusion' refers to the zonotope inclusion in \cite{Alamo2005a} for zonotopes and the CZ-inclusion (Theorem \ref{the:gzielim}) for constrained zonotopes. `Closest point' and `Change center' correspond to Propositions \ref{propo:closest} and \ref{propo:center}, respectively, which are LPs. For the latter, the bounds $\tilde{\bm{\xi}}^\text{L}, \tilde{\bm{\xi}}^\text{U}$ are obtained using Algorithm 1 in \cite{Scott2016}. Note that the interval hull of zonotopes does not require the solution of LPs (see Remark 3 in \cite{Kuhn1998}). In addition, we consider that each LP is solved at least with the performance of the simplex method presented in \cite{Kelner2006}, which is $O(N_dN_c^3)$ with $N_d$ and $N_c$ the number of decision variables and constraints, respectively. Note that these numbers can be inferred for each respective LP directly from Table \ref{tab:complexity}. For a detailed derivation of the computational complexities in Tables \ref{tab:complexityestimators} and \ref{tab:complexity}, please see the supplementary material.

\begin{table}[!htb]
	\scriptsize
	\centering
	\caption{Computational complexity $O(\cdot)$ of basic operations.}
	\begin{tabular}{c c c} \hline
		Operation & Zonotopes & Constrained zonotopes\\ \hline
	    Linear mapping & $nn_gn_r$ & $nn_gn_r$\\
	    Minkowski sum & $n$ & $n$ \\		
	    Generalized intersection & -- & $nn_gn_r + n_rn_{g_r}$\\
		Interval hull & $nn_g$ & $nn_g(n_g+n_c)^3$ \\
		Set inclusion & $nn_g$ & $nn_sn_g+n_c(n_g+n_c)^3+nn_g^2n_c$ \\
		Closest point & -- & $(n+n_g)(n+n_g+n_c)^3$ \\
		Change center& -- & $n_g(n+n_g)^3 + n_g^2n_c$ \\
		Generator reduction & $n^2n_g+k_gnn_g$ & $(n+n_c)^2n_g+k_g(n+n_c)n_g$ \\
		Constraint elimination & -- & $k_c(n_g+n_c)^3 + k_cnn_g^2$ \\
		\hline		
	\end{tabular} \normalsize
    \label{tab:complexity}
\end{table}


\bibliography{paper_bibliography}



\section*{Supplementary material}

# Supplementary material[★]

Brenner S. Rego[a,∗], Guilherme V. Raffo[a,b], Joseph K. Scott[c], Davide M. Raimondo[d]

[a]*Graduate Program in Electrical Engineering, Federal University of Minas Gerais, Belo Horizonte, MG 31270-901, Brazil*
[b]*Department of Electronics Engineering, Federal University of Minas Gerais, Belo Horizonte, MG 31270-901, Brazil*
[c]*Department of Chemical and Biomolecular Engineering, Clemson University, Clemson, SC, USA*
[d]*Department of Electrical, Computer and Biomedical Engineering, University of Pavia, Italy*

**Abstract**

This supplement details the derivation of the computational complexities of basic operations on zonotopes and constrained zonotopes, and the set-valued state estimators presented in the main paper.

*Keywords:* Nonlinear state estimation, Set-based computing, Reachability analysis, Convex polytopes

## 1. Computational complexity details: basic operations

*1.1. Zonotopes*

*1.1.1. Linear image*

Let $Z = \{\mathbf{G}_z, \mathbf{c}_z\} \subset \mathbb{R}^n$, $\mathbf{R} \in \mathbb{R}^{n_r \times n}$, with $\mathbf{c}_z \in \mathbb{R}^n$, $\mathbf{G}_z \in \mathbb{R}^{n \times n_g}$, and consider the linear image $\mathbf{R}Z = \{\mathbf{R}\mathbf{G}_z, \mathbf{R}\mathbf{c}_z\}$. Then,

$$O(\mathbf{R}Z) = O(\mathbf{R}\mathbf{c}_z) + O(\mathbf{R}\mathbf{G}_z) = O(nn_r) + O(nn_g n_r) = O(nn_g n_r).$$

*1.1.2. Minkowski sum*

Let $Z = \{\mathbf{G}_z, \mathbf{c}_z\} \subset \mathbb{R}^n$, $W = \{\mathbf{G}_w, \mathbf{c}_w\} \subset \mathbb{R}^n$, with $\mathbf{c}_z \in \mathbb{R}^n$, $\mathbf{G}_z \in \mathbb{R}^{n \times n_g}$, $\mathbf{c}_w \in \mathbb{R}^n$, $\mathbf{G}_z \in \mathbb{R}^{n \times n_{g_w}}$, and consider the Minkowski sum

$$Z \oplus W = \left\{ \begin{bmatrix} \mathbf{G}_z & \mathbf{G}_w \end{bmatrix}, \mathbf{c}_z + \mathbf{c}_w \right\}.$$

Then,

$$O(Z \oplus W) = O(\mathbf{c}_z + \mathbf{c}_w) = O(n).$$

*1.1.3. Generator reduction*

Let $Z = \{\mathbf{G}, \mathbf{c}\} \subset \mathbb{R}^n$, with $\mathbf{c} \in \mathbb{R}^n$, $\mathbf{G} \in \mathbb{R}^{n \times n_g}$. Let $k_g \leq (n_g - n)$ denote the number of generators to be eliminated. The computational complexity of generator reduction of zonotopes is the same presented in [1, 2]:

$$O(\text{eliminate } k_g \text{ generators from } Z) = O(n^2 n_g + k_g n_g n).$$

---

∗Corresponding author
*Email addresses:* brennersr7@ufmg.br (Brenner S. Rego), raffo@ufmg.br (Guilherme V. Raffo), jks9@clemson.edu (Joseph K. Scott), davide.raimondo@unipv.it (Davide M. Raimondo)

*1.1.4. Interval hull*

Let $Z = \{\mathbf{G}, \mathbf{c}\} \subset \mathbb{R}^n$, with $\mathbf{c} \in \mathbb{R}^n$, $\mathbf{G} \in \mathbb{R}^{n \times n_g}$. Consider the interval hull of Z given by $[\mathbf{c} - \zeta, \mathbf{c} + \zeta]$, where [3]

$$\zeta_i = \sum_{j=1}^{n_g} |\mathbf{G}_{ij}|.$$

Then,

$$O(\text{interval hull of } Z) = O(\mathbf{c} - \zeta) + O(\mathbf{c} + \zeta) + nO\left(\sum_{j=1}^{n_g} |\mathbf{G}_{ij}|\right) = O(\mathbf{c} + \zeta) + nO\left(\sum_{j=1}^{n_g} |\mathbf{G}_{ij}|\right)$$
$$= O(n) + nO(n_g) = O(nn_g).$$

*1.1.5. Zonotope inclusion*

Let $Z = \mathbf{p} \oplus \mathbf{M} B_\infty^{n_g}$, with $\mathbf{p} \in \mathbb{R}^n$, $\mathbf{M} \in \mathbb{IR}^{n \times n_g}$. Consider the zonotope inclusion $\diamond(Z) = \mathbf{p} \oplus [\text{mid}(\mathbf{M}) \ \mathbf{P}] B_\infty^{n_g+n}$, where $P_{ii} = \sum_{j=1}^{n_g} \frac{1}{2} \text{diam}(M_{ij})$, $P_{ij} = 0$ for $i \neq j$ [4]. Then,

$$O(\diamond(Z)) = O(\text{mid}(\mathbf{M})) + O(\mathbf{P}) = O(\text{mid}(\mathbf{M})) + nO(P_{ii})$$
$$= O(nn_g) + nO(n_g) = O(nn_g).$$

## 1.2. Constrained zonotopes

*1.2.1. Linear image*

Let $Z = \{\mathbf{G}_z, \mathbf{c}_z, \mathbf{A}_z, \mathbf{b}_z\} \subset \mathbb{R}^n$, $\mathbf{R} \in \mathbb{R}^{n_r \times n}$, with $\mathbf{c}_z \in \mathbb{R}^n$, $\mathbf{G}_z \in \mathbb{R}^{n \times n_g}$, $\mathbf{A}_z \in \mathbb{R}^{n_c \times n_g}$, $\mathbf{b}_z \in \mathbb{R}^{n_c}$, and consider the linear image $\mathbf{R}Z = \{\mathbf{R}\mathbf{G}_z, \mathbf{R}\mathbf{c}_z, \mathbf{A}_z, \mathbf{b}_z\}$. Then,

$$O(\mathbf{R}Z) = O(\mathbf{R}\mathbf{c}_z) + O(\mathbf{R}\mathbf{G}_z) = O(nn_r) + O(nn_g n_r) = O(nn_g n_r).$$

*1.2.2. Minkowski sum*

Let $Z = \{\mathbf{G}_z, \mathbf{c}_z, \mathbf{A}_z, \mathbf{b}_z\} \subset \mathbb{R}^n$, $W = \{\mathbf{G}_w, \mathbf{c}_w, \mathbf{A}_w, \mathbf{b}_w\} \subset \mathbb{R}^n$, with $\mathbf{c}_z \in \mathbb{R}^n$, $\mathbf{G}_z \in \mathbb{R}^{n \times n_g}$, $\mathbf{A}_z \in \mathbb{R}^{n_c \times n_g}$, $\mathbf{b}_z \in \mathbb{R}^{n_c}$, $\mathbf{c}_w \in \mathbb{R}^n$, $\mathbf{G}_w \in \mathbb{R}^{n \times n_{g_w}}$, $\mathbf{A}_w \in \mathbb{R}^{n_{c_w} \times n_{g_w}}$, $\mathbf{b}_w \in \mathbb{R}^{n_{c_w}}$, and consider the Minkowski sum

$$Z \oplus W = \left\{ \begin{bmatrix} \mathbf{G}_z & \mathbf{G}_w \end{bmatrix}, \mathbf{c}_z + \mathbf{c}_w, \begin{bmatrix} \mathbf{A}_z & \mathbf{0} \\ \mathbf{0} & \mathbf{A}_w \end{bmatrix}, \begin{bmatrix} \mathbf{b}_z \\ \mathbf{b}_w \end{bmatrix} \right\}.$$

Then,

$$O(Z \oplus W) = O(\mathbf{c}_z + \mathbf{c}_w) = O(n).$$

*1.2.3. Generalized intersection*

Let $Z = \{\mathbf{G}_z, \mathbf{c}_z, \mathbf{A}_z, \mathbf{b}_z\} \subset \mathbb{R}^n$, $Y = \{\mathbf{G}_y, \mathbf{c}_y, \mathbf{A}_y, \mathbf{b}_y\}$, $\mathbf{R} \in \mathbb{R}^{n_r \times n}$, with $\mathbf{G}_z \in \mathbb{R}^{n \times n_g}$, $\mathbf{c}_z \in \mathbb{R}^n$, $\mathbf{A}_z \in \mathbb{R}^{n_c \times n_g}$, $\mathbf{b}_z \in \mathbb{R}^{n_c}$, $\mathbf{G}_y \in \mathbb{R}^{n_r \times n_{g_r}}$, $\mathbf{c}_y \in \mathbb{R}^{n_r}$, $\mathbf{A}_y \in \mathbb{R}^{n_{c_r} \times n_{g_r}}$, $\mathbf{b}_y \in \mathbb{R}^{n_{c_r}}$. Consider the generalized intersection

$$Z \cap_{\mathbf{R}} Y = \left\{ \begin{bmatrix} \mathbf{G}_z & \mathbf{0} \end{bmatrix}, \mathbf{c}_z, \begin{bmatrix} \mathbf{A}_z & \mathbf{0} \\ \mathbf{0} & \mathbf{A}_y \\ \mathbf{R}\mathbf{G}_z & -\mathbf{G}_y \end{bmatrix}, \begin{bmatrix} \mathbf{b}_z \\ \mathbf{b}_y \\ \mathbf{c}_y - \mathbf{R}\mathbf{c}_z \end{bmatrix} \right\}.$$

Then,

$$O(Z \cap_{\mathbf{R}} Y) = O(\mathbf{R}\mathbf{G}_z) + O(-\mathbf{G}_y) + O(\mathbf{c}_y - \mathbf{R}\mathbf{c}_z) = O(nn_g n_r) + O(n_r n_{g_r}) + O(n_r + nn_r)$$
$$= O(nn_g n_r + n_r n_{g_r}).$$

### 1.2.4. Generator reduction

Let $Z = \{\mathbf{G}, \mathbf{c}, \mathbf{A}, \mathbf{b}\} \subset \mathbb{R}^n$, with $\mathbf{c} \in \mathbb{R}^n$, $\mathbf{G} \in \mathbb{R}^{n \times n_g}$, $\mathbf{A} \in \mathbb{R}^{n_c \times n_g}$, $\mathbf{b} \in \mathbb{R}^{n_c}$. Let $k_g \leq n_g - (n + n_c)$ denote the number of generators to be eliminated. The computational complexity of generator reduction of constrained zonotopes is the same presented in [1, 2], with reduction operated over the lifted zonotope

$$Z^+ = \left\{ \begin{bmatrix} \mathbf{G} \\ \mathbf{A} \end{bmatrix}, \begin{bmatrix} \mathbf{c} \\ -\mathbf{b} \end{bmatrix} \right\}.$$

Therefore,

$$O(\text{eliminate } k_g \text{ generators from } Z) = O((n + n_c)^2 n_g + k_g n_g (n + n_c)).$$

### 1.2.5. Constraint elimination

Let $Z = \{\mathbf{G}, \mathbf{c}, \mathbf{A}, \mathbf{b}\} \subset \mathbb{R}^n$, with $\mathbf{c} \in \mathbb{R}^n$, $\mathbf{G} \in \mathbb{R}^{n \times n_g}$, $\mathbf{A} \in \mathbb{R}^{n_c \times n_g}$, $\mathbf{b} \in \mathbb{R}^{n_c}$. Let $k_c \leq n_c$ denote the number of constraints to be eliminated from $Z$. Following the constraint elimination algorithm in [1], for each eliminated constraint the remaining constraints are first preconditioned through Gauss-Jordan elimination with full pivoting and then subjected to a rescaling procedure before the next constraint is eliminated. Therefore,

$$O(\text{eliminate 1 constraint from } Z) = O(\text{pre-conditioning}) + O(\text{rescaling}) + O(\text{eliminate the constraint}).$$

From the Appendix in [1], we have that

$$O(\text{pre-conditioning}) + O(\text{rescaling}) = O(n_c^2 n_g + n_c n_g^2).$$

To eliminate the constraint, one of the generators of $Z$ is first chosen by minimizing an approximation of the Hausdorff distance between $Z$ and the set obtained by constraint elimination. This requires solving (A.8) in [1] for each $j \in \{1, 2, \ldots, n_g\}$, which has the reported complexity $O((n_g + n_c)^3)$, and computing the Hausdorff distance $H_j^* = \|\mathbf{G}\mathbf{d}_j^*\|_2^2 + \|\mathbf{d}_j^*\|_2^2$, with $\mathbf{d}_j^*$ the solution of (A.8) for each $j$. This is done with the following complexity:

$$O(\text{eliminate the constraint}) = O(\text{Hausdorff distance minimization}) + n_g O(H_j^*)$$
$$= O((n_g + n_c)^3) + n_g O(n n_g) = O((n_g + n_c)^3 + n n_g^2).$$

Then,

$$O(\text{eliminate 1 constraint from } Z) = O(n_c^2 n_g + n_c n_g^2) + O((n_g + n_c)^3 + n n_g^2) = O((n_g + n_c)^3 + n n_g^2).$$

Finally, the complexity of eliminating $k_c$ constraints from $Z$ can be bounded by

$$O(\text{eliminate } k_c \text{ constraints from } Z) = k_c O(\text{eliminate 1 constraint from } Z) = O(k_c (n_g + n_c)^3 + k_c n n_g^2).$$

### 1.2.6. Interval hull (Property 1)

Let $Z = \{\mathbf{G}, \mathbf{c}, \mathbf{A}, \mathbf{b}\} \subset \mathbb{R}^n$, with $\mathbf{c} \in \mathbb{R}^n$, $\mathbf{G} \in \mathbb{R}^{n \times n_g}$, $\mathbf{A} \in \mathbb{R}^{n_c \times n_g}$, $\mathbf{b} \in \mathbb{R}^{n_c}$. Consider the interval hull of $Z$ computed as in Property 1. In the standard form, each LP in Property 1 has $N_d = n_g + 1$ decision variables and $N_c = 2n_g + n_c + 1$ constraints. Then,

$$O(\text{interval hull of } Z) = 2n O(N_d N_c^3) = O(n N_d N_c^3) = O(n(1 + n_g)(1 + 2n_g + n_c)^3) = O(n n_g (n_g + n_c)^3).$$

### 1.2.7. CZ-inclusion (Theorem 1)

Let $X = \mathbf{p} \oplus \mathbf{M} B_\infty(\mathbf{A}, \mathbf{b}) \subset \mathbb{R}^m$, $\mathbf{J} \in \mathbb{IR}^{n \times m}$, $\bar{X} = \bar{\mathbf{p}} \oplus \bar{\mathbf{M}} B_\infty^{\bar{n}_g} \supseteq X$, $\mathbf{m} \supseteq (\mathbf{J} - \text{mid}(\mathbf{J}))\bar{\mathbf{p}}$, $\triangleleft(\mathbf{J}, X) = \text{mid}(\mathbf{J}) X \oplus \mathbf{P} B_\infty^n$, with $\mathbf{p} \in \mathbb{R}^m$, $\mathbf{M} \in \mathbb{R}^{m \times n_g}$, $\mathbf{A} \in \mathbb{R}^{n_c \times n_g}$, $\mathbf{b} \in \mathbb{R}^{n_c}$, $P_{ii} = \frac{1}{2}\text{diam}(m_i) + \frac{1}{2} \sum_{j=1}^{\bar{n}_g} \sum_{k=1}^{m} \text{diam}(J_{ik})|\bar{M}_{kj}|$, $P_{ij} = 0$ for $i \neq j$. Assume that $\bar{X}$ is obtained by eliminating $n_c$ constraints from $X$, and $\mathbf{m}$ is computed using interval arithmetic. Then,

$$O(\triangleleft(\mathbf{J}, X)) = O(\text{mid}(\mathbf{J}) X \oplus \mathbf{P} B_\infty^n) + O((\mathbf{J} - \text{mid}(\mathbf{J}))\bar{\mathbf{p}}) + n O\left( (1/2)\text{diam}(m_i) + \frac{1}{2} \sum_{j=1}^{\bar{n}_g} \sum_{k=1}^{m} \text{diam}(J_{ik})|\bar{M}_{kj}| \right)$$

$$+ O(\text{eliminate } n_c \text{ constraints from } X)$$
$$= O(nm + nmn_g + n) + O(nm + nm) + nO(\bar{n}_g m) + O(n_c(n_g + n_c)^3 + n_c m n_g^2)$$
$$= O(nmn_g) + O(nm) + O(nmn_g) + O(n_c(n_g + n_c)^3 + mn_c n_g^2)$$
$$= O(nmn_g + n_c(n_g + n_c)^3 + mn_g^2 n_c).$$

### 1.2.8. Closest point (Proposition 1)

Let $Z = \{\mathbf{G}, \mathbf{c}, \mathbf{A}, \mathbf{b}\} \subset \mathbb{R}^n$, $\mathbf{h} \in \mathbb{R}^n$, with $\mathbf{c} \in \mathbb{R}^n$, $\mathbf{G} \in \mathbb{R}^{n \times n_g}$, $\mathbf{A} \in \mathbb{R}^{n_c \times n_g}$, $\mathbf{b} \in \mathbb{R}^{n_c}$. In the standard form, the LP in Proposition 1 has $N_d = n + n_g + 1$ decision variables and $N_c = 2n + 2n_g + n_c + 1$ constraints. Then,

$$O(\text{Proposition 1}) = O(\mathbf{c} + \mathbf{G}\boldsymbol{\xi}^*) + O(N_d N_c^3) = O(n + nn_g) + O((n + n_g + 1)(2n + 2n_g + n_c + 1)^3)$$
$$= O((n + n_g)(n + n_g + n_c)^3).$$

### 1.2.9. Change center (Proposition 2)

Let
$$Z = \{\mathbf{G}, \mathbf{c}, \mathbf{A}, \mathbf{b}\} = \left\{ \begin{bmatrix} \mathbf{G}\mathbf{E}_r & \mathbf{0} \end{bmatrix}, \mathbf{h}, \begin{bmatrix} \mathbf{A}\mathbf{E}_r & \mathbf{0} \\ \mathbf{0} & \mathbf{A} \\ \mathbf{G}\mathbf{E}_r & -\mathbf{G} \end{bmatrix}, \begin{bmatrix} \mathbf{b} - \mathbf{A}\boldsymbol{\xi}_m \\ \mathbf{b} \\ \mathbf{c} - \mathbf{h} \end{bmatrix} \right\} \subset \mathbb{R}^n,$$

with $\mathbf{c} \in \mathbb{R}^n$, $\mathbf{G} \in \mathbb{R}^{n \times n_g}$, $\mathbf{A} \in \mathbb{R}^{n_c \times n_g}$, $\mathbf{b} \in \mathbb{R}^{n_c}$, $\mathbf{h} \in \mathbb{R}^n$, $\mathbf{E}_r = \frac{1}{2}\text{diag}(\boldsymbol{\xi}^U - \boldsymbol{\xi}^L)$, $\boldsymbol{\xi}_m = \frac{1}{2}(\boldsymbol{\xi}^U + \boldsymbol{\xi}^L)$. Consider $(\boldsymbol{\xi}^L, \boldsymbol{\xi}^U)$ computed by solving the LP in Proposition 2, in which $(\tilde{\boldsymbol{\xi}}^L, \tilde{\boldsymbol{\xi}}^U)$ is obtained by Algorithm 1 in the Appendix in [1]. In the standard form, the LP in Proposition 2 has $N_d = 3n_g$ decision variables and $N_c = n + 4n_g$ constraints. Then,

$$O(\text{Proposition 2}) = O(\mathbf{G}\mathbf{E}_r) + O(\mathbf{A}\mathbf{E}_r) + O(-\mathbf{G}) + O(\mathbf{b} - \mathbf{A}\boldsymbol{\xi}_m) + O(\mathbf{c} - \mathbf{h})$$
$$+ O((1/2)\text{diag}(\boldsymbol{\xi}^U - \boldsymbol{\xi}^L)) + O((1/2)(\boldsymbol{\xi}^U + \boldsymbol{\xi}^L)) + O(N_d N_c^3) + O(\text{Algorithm 1})$$
$$= O(nn_g^2) + O(n_c n_g^2) + O(nn_g) + O(n_c) + O(n_c n_g) + O(n) + O(n_g) + O(n_g)$$
$$+ O(N_d N_c^3) + O(\text{Algorithm 1})$$
$$= O(nn_g^2) + O(n_c n_g^2) + O(N_d N_c^3) + O(\text{Algorithm 1})$$
$$= O(nn_g^2) + O(n_c n_g^2) + O((3n_g)(n + 4n_g)^3) + O(\text{Algorithm 1})$$
$$= O(nn_g^2) + O(n_c n_g^2) + O(n_g(n + n_g)^3) + O(\text{Algorithm 1})$$
$$= O(n_g(n + n_g)^3 + n_g^2 n_c) + O(n_g^2 n_c) = O(n_g(n + n_g)^3 + n_g^2 n_c).$$

## 2. Computational complexity details: ZMV

### 2.1. Prediction step

Let $\boldsymbol{\mu} : \mathbb{R}^n \times \mathbb{R}^{n_w} \to \mathbb{R}^n$ be continuously differentiable and let $\nabla_x \boldsymbol{\mu}$ denote the gradient of $\boldsymbol{\mu}$ with respect to its first argument. Let $\hat{X}_{k-1} = \{\mathbf{G}_x, \mathbf{c}_x\} \subset \mathbb{R}^n$, and $W = \{\mathbf{G}_w, \mathbf{c}_w\} \subset \mathbb{R}^{n_w}$ with $\mathbf{c}_x \in \mathbb{R}^n$, $\mathbf{G}_x \in \mathbb{R}^{n \times n_g}$, $\mathbf{c}_w \in \mathbb{R}^{n_w}$, and $\mathbf{G}_w \in \mathbb{R}^{n_w \times n_{g_w}}$. Let $Z$ be a zonotope such that $\boldsymbol{\mu}(\mathbf{c}_x, W) \subseteq Z$ and let $\mathbf{M} \in \mathbb{IR}^{n \times n_g}$ be an interval matrix satisfying $\nabla_x^T \boldsymbol{\mu}(\hat{X}_{k-1}, W)\mathbf{G}_\mathbf{x} \subseteq \mathbf{M}$. According to Theorem 4 in [4],

$$\boldsymbol{\mu}(\hat{X}_{k-1}, W) \subseteq \bar{X}_k = Z \oplus \diamond(\mathbf{M} B_\infty^{n_g}),$$

in which $\bar{X}_k$ has $n_g + n_{g_w} + 2n$ generators if $Z$ is also computed by the mean value extension. Therefore,

$$O(\text{ZMV}_{\text{prediction}}) = O(Z \oplus \diamond(\mathbf{M} B_\infty^{n_g})) + O(\diamond(\mathbf{M} B_\infty^{n_g})) + O(Z) + O(\text{interval hull of } \hat{X}_{k-1}, W)$$
$$= O(n) + O(nn_g + n^2 n_g) + O(n_w n_{g_w} + nn_w n_{g_w}) + O(nn_g + n_w n_{g_w})$$
$$= O(n^2 n_g + nn_w n_{g_w}).$$

## 2.2. Update step

Let $\mathbf{y}_k = \mathbf{C}\mathbf{x}_k + \mathbf{D}_u\mathbf{u}_k + \mathbf{D}_v\mathbf{v}_k$, with $\mathbf{y}_k \in \mathbb{R}^{n_y}$, $\mathbf{x}_k \in \mathbb{R}^n$, $\mathbf{u}_k \in \mathbb{R}^{n_u}$, $\mathbf{v}_k \in \mathbb{R}^{n_v}$. Consider the zonotopes $\bar{X}_k = \{\mathbf{G}_x, \mathbf{c}_x\} \subset \mathbb{R}^n$, $V = \{\mathbf{G}_v, \mathbf{c}_v\} \subset \mathbb{R}^{n_v}$ with $\mathbf{c}_x \in \mathbb{R}^n$, $\mathbf{G}_x \in \mathbb{R}^{n \times \bar{n}_g}$, $\mathbf{c}_v \in \mathbb{R}^{n_v}$, and $\mathbf{G}_v \in \mathbb{R}^{n_v \times n_{g_v}}$. For simplicity, define $m_g = n_g + n_{g_w}$. Following the intersection method proposed in [5], the update step is computed iteratively using one row of $\mathbf{y}_k$ at a time, as described below.

*Algorithm: Update step.*

(1) Assign $\tilde{X}_k = \{\tilde{\mathbf{G}}_x, \tilde{\mathbf{c}}_x\} \leftarrow \bar{X}_k$, $i \leftarrow 1$.
(2) Assign $\tilde{y} \leftarrow$ *i*-th row of $\mathbf{y}_k$, $\mathbf{p}^T \leftarrow$ *i*-th row of $\mathbf{C}$, $\mathbf{d}_u^T \leftarrow$ *i*-th row of $\mathbf{D}_u$, $\mathbf{d}_v^T \leftarrow$ *i*-th row of $\mathbf{D}_v$.
(3) Compute the strip $S \leftarrow \{\mathbf{x} : |\mathbf{p}^T\mathbf{x} - d| \leq \sigma\}$, where $d = \tilde{y} - \mathbf{d}_u^T\mathbf{u}_k - \mathbf{d}_v\mathbf{c}_v$, $\sigma = \sum_{\ell=1}^{n_{g_v}} \sum_{j=1}^{n_v} |(\mathbf{D}_v)_{ij}(\mathbf{G}_v)_{j\ell}|$.
(4) Compute the tight strip $\tilde{S} \leftarrow \{\mathbf{x} : |\mathbf{p}^T\mathbf{x} - \tilde{d}| \leq \tilde{\sigma}\}$, where

$$\tilde{\sigma} = \frac{1}{2}(\bar{\sigma}^U - \bar{\sigma}^L), \quad \tilde{d} = d + \frac{1}{2}(\bar{\sigma}^U + \bar{\sigma}^L),$$

with $\bar{\sigma}^L = \max\{-\sigma, \mathbf{p}^T\tilde{\mathbf{c}}_x - \|\tilde{\mathbf{G}}_x^T\mathbf{p}\|_1 - d\}$, $\bar{\sigma}^U = \min\{\sigma, \mathbf{p}^T\tilde{\mathbf{c}}_x + \|\tilde{\mathbf{G}}_x^T\mathbf{p}\|_1 - d\}$.

(5) Assign $\tilde{X}_k \leftarrow \tilde{X}_k \cap \tilde{S} = \{G(j^*), c(j^*)\}$, where

$$j^* = \arg\min_j \det(\mathbf{G}(j)\mathbf{G}(j)^T),$$

with $j \in \{0, 1, \ldots, \bar{n}_g\}$,

$$\mathbf{c}(j) = \begin{cases} \tilde{\mathbf{c}}_x + \left(\dfrac{\tilde{d} - \mathbf{p}^T\tilde{\mathbf{c}}_x}{\mathbf{p}^T(\tilde{\mathbf{G}}_x)_{:,j}}\right)(\tilde{\mathbf{G}}_x)_{:,j} & \text{if } 1 \leq j \leq \bar{n}_g \text{ and } \mathbf{p}^T(\tilde{\mathbf{G}}_x)_{:,j} \neq 0 \\ \tilde{\mathbf{c}}_x & \text{otherwise} \end{cases}$$

$$\mathbf{G}(j) = \begin{cases} [\mathbf{g}_1^j \; \mathbf{g}_2^j \; \ldots \; \mathbf{g}_{\bar{n}_g}^j] & \text{if } 1 \leq j \leq \bar{n}_g \text{ and } \mathbf{p}^T(\tilde{\mathbf{G}}_x)_{:,j} \neq 0 \\ \tilde{\mathbf{G}}_x & \text{otherwise} \end{cases}$$

$$\mathbf{g}_\ell^j = \begin{cases} (\tilde{\mathbf{G}}_x)_{:,\ell} - \left(\dfrac{\mathbf{p}^T(\tilde{\mathbf{G}}_x)_{:,\ell}}{\mathbf{p}^T(\tilde{\mathbf{G}}_x)_{:,j}}\right)(\tilde{\mathbf{G}}_x)_{:,j} & \text{if } \ell \neq j \\ \left(\dfrac{\tilde{\sigma}}{\mathbf{p}^T(\tilde{\mathbf{G}}_x)_{:,j}}\right)(\tilde{\mathbf{G}}_x)_{:,j} & \text{if } \ell = j \end{cases}$$

(6) If $i < n_y$, go to Step 2 and assign $i \leftarrow i + 1$. Otherwise, assign $\hat{X}_k \leftarrow \tilde{X}_k$.

The computational complexity is

$$O(\text{ZMV}_{\text{update}}) = n_y O(\text{update with the } i\text{-th row of } \mathbf{y}_k)$$
$$= n_y \left(O(\text{Step (3)}) + O(\text{Step (4)}) + O(\text{Step (5)})\right).$$

But,

$$O(\text{Step (3)}) = O(1) + O(n_u) + O(n_v) + O(n_v n_{g_v}) = O(n_u + n_v n_{g_v}),$$

$$O(\text{Step (4)}) = O(1) + O(n) + O(n\bar{n}_g) = O(n\bar{n}_g).$$

$$O(\text{Step (5)}) = \bar{n}_g O(n^3) + \bar{n}_g O(n^2 \bar{n}_g) + \bar{n}_g O(n) + \bar{n}_g O(n\bar{n}_g) = O(n^3 \bar{n}_g + n^2 \bar{n}_g^2).$$

Since $\bar{n}_g = m_g + 2n$, hence

$$O(\text{ZMV}_{\text{update}}) = n_y \left(O(n_u + n_v n_{g_v}) + O(n\bar{n}_g) + O(n^3 \bar{n}_g + n^2 \bar{n}_g^2 + n\bar{n}_g^2)\right)$$
$$= n_y \left(O(n^3 \bar{n}_g + n^2 \bar{n}_g^2 + n_u + n_v n_{g_v})\right)$$
$$= O(n_y(n^3(m_g + n) + n^2(m_g + n)^2 + n_u + n_v n_{g_v})).$$

5

*2.3. Order reduction*

Let $n_g$ denote the desired number of generators of $\hat{X}_k$ after order reduction, i.e., $\hat{X}_k$ must have the same complexity as $\hat{X}_{k-1}$. Order reduction is performed by eliminating $k_g$ generators from $\hat{X}_k$. For simplicity, define $m_g = n_g + n_{g_w}$. The prediction and update steps of the ZMV lead to $\hat{X}_k$ with $m_g + 2n$ generators. Therefore, $k_g = n_{g_w} + 2n$. Consequently,

$$O(\text{ZMV}_{\text{reduction}}) = O(\text{eliminate } k_g \text{ generators from } \hat{X}_k)$$
$$= O(n^2(m_g + n) + k_g(m_g + n)n)$$
$$= O(n^2(m_g + n) + (n_{g_w} + n)(m_g + n)n) = O(n^2(m_g + n) + n(n_{g_w} + n)(m_g + n)).$$

## 3. Computational complexity details: ZFO

Let $\boldsymbol{\eta} : \mathbb{R}^m \to \mathbb{R}^n$ be of class $C^2$ and let $\mathbf{z} = (\mathbf{x}_{k-1}, \mathbf{w}_{k-1}) \in \mathbb{R}^m$ denote its argument, where $m = n + n_w$. Let $Z = \hat{X}_{k-1} \times W = \{\mathbf{G}, \mathbf{c}\} \subset \mathbb{R}^m$ be a zonotope with $m_g = n_g + n_{g_w}$ generators. For each $q = 1, 2, \ldots, n$, let $\mathbf{Q}^{[q]} \in \mathbb{IR}^{m \times m}$ and $\tilde{\mathbf{Q}}^{[q]} \in \mathbb{IR}^{m_g \times m_g}$ be interval matrices satisfying $\mathbf{Q}^{[q]} \supseteq \mathbf{H}\eta_q(Z)$ and $\tilde{\mathbf{Q}}^{[q]} \supseteq \mathbf{G}^T \mathbf{Q}^{[q]} \mathbf{G}$. Moreover, define $\tilde{\mathbf{c}}$, $\tilde{\mathbf{G}}$, $\tilde{\mathbf{G}}_\mathbf{d}$ as in Theorem 3. Then, the method in [6] ensures that

$$\boldsymbol{\eta}(Z) \subseteq \bar{X}_k = \boldsymbol{\eta}(\mathbf{c}) \oplus \nabla^T \boldsymbol{\eta}(\mathbf{c})(\mathbf{G}B_\infty^{m_g}) \oplus \tilde{\mathbf{c}} \oplus \begin{bmatrix} \tilde{\mathbf{G}} & \tilde{\mathbf{G}}_\mathbf{d} \end{bmatrix} B_\infty^{\tilde{m}_g},$$

where $\tilde{m}_g = (1/2)m_g^2 + (3/2)m_g + n$, and therefore $\bar{X}_k$ has $(1/2)m_g^2 + (5/2)m_g + n$ generators. Thus,

$$O(\text{CZFO}_{\text{prediction}}) = O(\boldsymbol{\eta}(\mathbf{c}) \oplus \nabla^T \boldsymbol{\eta}(\mathbf{c})(\mathbf{G}B_\infty^{m_g}) \oplus \tilde{\mathbf{c}} \oplus \begin{bmatrix} \tilde{\mathbf{G}} & \tilde{\mathbf{G}}_\mathbf{d} \end{bmatrix} B_\infty^{\tilde{m}_g})$$
$$= O(n) + O(\tilde{\mathbf{c}}) + O(\tilde{\mathbf{G}}) + O(\tilde{\mathbf{G}}_\mathbf{d}) + nO(\tilde{\mathbf{Q}}^{[q]}) + O(\text{interval hull of } Z)$$
$$= O(n) + O(nm_g) + O(nm_g^2) + O(nm_g^2) + nO(m^2 m_g + m m_g^2) + O(m m_g)$$
$$= O(n(m^2 m_g + m m_g^2)).$$

*3.1. Update step*

Let $\mathbf{y}_k = \mathbf{C}\mathbf{x}_k + \mathbf{D}_u \mathbf{u}_k + \mathbf{D}_v \mathbf{v}_k$, with $\mathbf{y}_k \in \mathbb{R}^{n_y}$, $\mathbf{x}_k \in \mathbb{R}^n$, $\mathbf{u}_k \in \mathbb{R}^{n_u}$, and $\mathbf{v}_k \in \mathbb{R}^{n_v}$. Consider the zonotopes $\bar{X}_k = \{\mathbf{G}_x, \mathbf{c}_x\} \subset \mathbb{R}^n$ and $V = \{\mathbf{G}_v, \mathbf{c}_v\} \subset \mathbb{R}^{n_v}$ with $\mathbf{c}_x \in \mathbb{R}^n$, $\mathbf{G}_x \in \mathbb{R}^{n \times \bar{n}_g}$, $\mathbf{c}_v \in \mathbb{R}^{n_v}$, and $\mathbf{G}_v \in \mathbb{R}^{n_v \times n_{g_v}}$. For simplicity, define $m_g = n_g + n_{g_w}$. The ZFO also follows the intersection method proposed in [5]. Therefore, the update step is equivalent to the ZMV update, but with $\bar{n}_g = (1/2)m_g^2 + (5/2)m_g + n$. Hence

$$O(\text{ZFO}_{\text{update}}) = n_y \left( O(n^3 \bar{n}_g + n^2 \bar{n}_g^2 + n_u + n_v n_{g_v}) \right)$$
$$= O(n_y(n^3(m_g^2 + n) + n^2(m_g^2 + n)^2 + n_u + n_v n_{g_v})).$$

*3.2. Order reduction*

Let $n_g$ denote the desired number of generators of $\hat{X}_k$ after order reduction, i.e., $\hat{X}_k$ must have the same complexity as $\hat{X}_{k-1}$. Order reduction is performed by eliminating $k_g$ generators from $\hat{X}_k$. For simplicity, define $m_g = n_g + n_{g_w}$. The prediction and update steps of ZFO lead to $\hat{X}_k$ with $(1/2)m_g^2 + (5/2)m_g + n$ generators. Therefore, $k_g = (1/2)m_g^2 + (5/2)m_g + n - n_g$. Consequently,

$$O(\text{ZFO}_{\text{reduction}}) = O(\text{eliminate } k_g \text{ generators from } \hat{X}_k)$$
$$= O(n^2(m_g^2 + m_g + n) + k_g(m_g^2 + m_g + n)n)$$
$$= O(n^2(m_g^2 + n) + (m_g^2 + n)(m_g^2 + n)n) = O(n^2(m_g^2 + n) + n(m_g^2 + n)^2).$$

# 4. Computational complexity details: CZMV

## 4.1. Prediction step (Theorem 2)

Let $\mu : \mathbb{R}^n \times \mathbb{R}^{n_w} \to \mathbb{R}^n$ be continuously differentiable and let $\nabla_x \mu$ denote the gradient of $\mu$ with respect to its first argument. Let $\hat{X}_{k-1} = \{\mathbf{G}_x, \mathbf{c}_x, \mathbf{A}_x, \mathbf{b}_x\} \subset \mathbb{R}^n$ and $W = \{\mathbf{G}_w, \mathbf{c}_w, \mathbf{A}_w, \mathbf{b}_w\} \subset \mathbb{R}^{n_w}$ with $\mathbf{c}_x \in \mathbb{R}^n$, $\mathbf{G}_x \in \mathbb{R}^{n \times n_g}$, $\mathbf{A}_x \in \mathbb{R}^{n_c \times n_g}$, $\mathbf{b}_x \in \mathbb{R}^{n_c}$, $\mathbf{c}_w \in \mathbb{R}^{n_w}$, $\mathbf{G}_w \in \mathbb{R}^{n_w \times n_{g_w}}$, $\mathbf{A}_w \in \mathbb{R}^{n_{c_w} \times n_{g_w}}$, and $\mathbf{b}_w \in \mathbb{R}^{n_{c_w}}$, and choose any $\mathbf{h} \in \hat{X}_{k-1}$. Let $Z$ be a constrained zonotope such that $\mu(\mathbf{h}, W) \subseteq Z$ and let $\mathbf{J} \in \mathbb{IR}^{n \times n}$ be an interval matrix satisfying $\nabla_x^T \mu(\hat{X}_{k-1}, W) \subseteq \mathbf{J}$. Theorem 2 ensures that

$$\mu(\hat{X}_{k-1}, W) \subseteq \bar{X}_k = Z \oplus \triangleleft(\mathbf{J}, \hat{X}_{k-1} - \mathbf{h}),$$

where $\bar{X}_k$ has $n_g + n_{g_w} + 2n$ generators and $n_c + n_{c_w}$ constraints if $Z$ is also computed by the mean value extension. Therefore,

$$\begin{aligned}
O(\text{CZMV}_{\text{prediction}}) &= O(Z \oplus \triangleleft(\mathbf{J}, \hat{X}_{k-1} - \mathbf{h})) + O(\triangleleft(\mathbf{J}, \hat{X}_{k-1} - \mathbf{h})) + O(Z) + O(\text{interval hull of } \hat{X}_{k-1}, W) \\
&= O(n) + O(n^2 n_g + n_c(n_g + n_c)^3 + n n_g^2 n_c) + O(n n_w n_{g_w} + n_{c_w}(n_{g_w} + n_{c_w})^3 + n_w n_{g_w}^2 n_{c_w}) \\
&\quad + O(n n_g (n_g + n_c)^3 + n_w n_{g_w}(n_{g_w} + n_{c_w})^3) \\
&= O(n^2 n_g + n n_w n_{g_w} + (n n_g + n_c)(n_g + n_c)^3 + (n_w n_{g_w} + n_{c_w})(n_{g_w} + n_{c_w})^3).
\end{aligned}$$

## 4.2. Update step

Let $\bar{X}_k = \{\mathbf{G}_x, \mathbf{c}_x, \mathbf{A}_x, \mathbf{b}_x\} \subset \mathbb{R}^n$, $V = \{\mathbf{G}_v, \mathbf{c}_v, \mathbf{A}_v, \mathbf{b}_v\} \subset \mathbb{R}^{n_w}$, $\mathbf{u}_k \in \mathbb{R}^{n_u}$, and $\mathbf{y}_k \in \mathbb{R}^{n_y}$ with $\mathbf{c}_x \in \mathbb{R}^n$, $\mathbf{G}_x \in \mathbb{R}^{n \times \bar{n}_g}$, $\mathbf{A}_x \in \mathbb{R}^{\bar{n}_c \times \bar{n}_g}$, $\mathbf{b}_x \in \mathbb{R}^{\bar{n}_c}$, $\mathbf{c}_v \in \mathbb{R}^{n_v}$, $\mathbf{G}_v \in \mathbb{R}^{n_v \times n_{g_v}}$, $\mathbf{A}_v \in \mathbb{R}^{n_{c_v} \times n_{g_v}}$, and $\mathbf{b}_w \in \mathbb{R}^{n_{c_v}}$. The update step is performed according to

$$\hat{X}_k = \bar{X}_k \cap_{\mathbf{C}} ((\mathbf{y}_k - \mathbf{D}_u \mathbf{u}_k) \oplus (-\mathbf{D}_v V)),$$

where $\mathbf{C} \in \mathbb{R}^{n_y \times n}$, $\mathbf{D}_u \in \mathbb{R}^{n_y \times n_u}$, and $\mathbf{D}_v \in \mathbb{R}^{n_y \times n_v}$. As a result, from the prediction step of CZMV, $\bar{X}_k$ has $\bar{n}_g = n_g + n_{g_w} + 2n$ generators and $\bar{n}_c = n_c + n_{c_w}$ constraints. Therefore,

$$\begin{aligned}
O(\text{CZMV}_{\text{update}}) &= O(\bar{X}_k \cap_{\mathbf{C}} ((\mathbf{y}_k - \mathbf{D}_u \mathbf{u}_k) \oplus (-\mathbf{D}_v V))) \\
&= O(n(n_g + n_{g_w} + 2n)n_y + n_y n_{g_v}) + O(n_y) + O(n_y + n_y n_u) + O(n_y n_v n_{g_v}) \\
&= O(n_y n(m_g + n) + n_y n_u + n_y n_v n_{g_v}),
\end{aligned}$$

where $m_g = n_g + n_{g_w}$.

## 4.3. Order reduction

Let $n_g$ and $n_c$ denote the desired number of generators and constraints of $\hat{X}_k$ after order reduction, i.e., $\hat{X}_k$ must have the same complexity as $\hat{X}_{k-1}$. Order reduction is performed by eliminating $k_c$ constraints and then $k_g$ generators from $\hat{X}_k$. For simplicity, define $m_g = n_g + n_{g_w}$, $m_c = n_c + n_{c_w}$, $\delta_n = n - n_y$, $\delta_w = n_{g_w} - n_{c_w}$, and $\delta_v = n_{g_v} - n_{c_v}$. The prediction and update steps of CZMV lead to $\hat{X}_k$ with $n_g + n_{g_w} + 2n + n_{g_v}$ generators and $n_c + n_{c_w} + n_{c_v} + n_y$ constraints. Therefore, $k_c = n_{c_w} + n_{c_v} + n_y$, and $k_g = n_{g_w} + 2n + n_{g_v} - k_c = n_{g_w} + 2n + n_{g_v} - n_{c_w} - n_{c_v} - n_y$. Consequently,

$$\begin{aligned}
O(\text{eliminate } k_c \text{ constraints from } \hat{X}_k) &= O(k_c(m_g + 2n + n_{g_v} + m_c + n_{c_v} + n_y)^3 \\
&\quad + k_c n(m_g + 2n + n_{g_v})^2) \\
&= O((n_{c_w} + n_{c_v} + n_y)(m_g + m_c + n_{g_v} + n_{c_v} + n + n_y)^3 \\
&\quad + (n_{c_w} + n_{c_v} + n_y)n(m_g + n_{g_v} + n)^2) \\
&= O((n_{c_w} + n_{c_v} + n_y)(m_g + m_c + n_{g_v} + n_{c_v} + n + n_y)^3.
\end{aligned}$$

Note that by eliminating $k_c$ constraints from $\hat{X}_k$, the same number of generators (i.e., $k_c$) is also removed [1]. Hence,

$$\begin{aligned}
O(\text{elim. } k_g \text{ gen. from } \hat{X}_k \text{ after const. elim.}) &= O((n + n_c)^2(m_g + 2n + n_{g_v} - k_c) \\
&\quad + k_g(m_g + 2n + n_{g_v} - k_c)(n + n_c))
\end{aligned}$$

7

$$= O((n + n_c)^2(\delta_n + \delta_w + \delta_v + n_g)$$
$$+ (\delta_n + \delta_w + \delta_v)(\delta_n + \delta_w + \delta_v + n_g)(n + n_c)).$$

Finally,

$$O(\text{CZMV}_{\text{reduction}}) = O(\text{eliminate } k_c \text{ constraints from } \hat{X}_k) + O(\text{elim. } k_g \text{ gen. from } \hat{X}_k \text{ after const. elim.})$$
$$= O((n_{c_w} + n_{c_v} + n_y)(m_g + m_c + n_{g_v} + n_{c_v} + n + n_y)^3 + (n + n_c)^2(n_g + \delta_n + \delta_w + \delta_v)$$
$$+ (n + n_c)(\delta_n + \delta_w + \delta_v)(n_g + \delta_n + \delta_w + \delta_v)).$$

## 5. Computational complexity details: CZFO

### 5.1. Prediction step

Let $\eta : \mathbb{R}^m \to \mathbb{R}^n$ be of class $C^2$ and let $\mathbf{z} = (\mathbf{x}_{k-1}, \mathbf{w}_{k-1}) \in \mathbb{R}^m$ denote its argument, where $m = n + n_w$. Let $Z = \hat{X}_{k-1} \times W = \{\mathbf{G}, \mathbf{c}, \mathbf{A}, \mathbf{b}\} \subset \mathbb{R}^m$ be a constrained zonotope with $m_g = n_g + n_{g_w}$ generators and $m_c = n_c + n_{c_w}$ constraints. For each $q = 1, 2, \ldots, n$, let $\mathbf{Q}^{[q]} \in \mathbb{IR}^{m \times m}$ and $\tilde{\mathbf{Q}}^{[q]} \in \mathbb{IR}^{m_g \times m_g}$ be interval matrices satisfying $\mathbf{Q}^{[q]} \supseteq \mathbf{H}\eta_q(Z)$ and $\tilde{\mathbf{Q}}^{[q]} \supseteq \mathbf{G}^T \mathbf{Q}^{[q]} \mathbf{G}$. Moreover, define $\tilde{\mathbf{c}}, \tilde{\mathbf{G}}, \tilde{\mathbf{G}}_{\mathbf{d}}, \tilde{\mathbf{A}}, \tilde{\mathbf{b}}$, and $\mathbf{L}$ as in Theorem 3. Choosing any $\mathbf{h} \in Z$, Theorem 3 ensures that

$$\eta(Z) \subseteq \bar{X}_k = \eta(\mathbf{h}) \oplus \nabla^T \eta(\mathbf{h})(Z - \mathbf{h}) \oplus R,$$

where $R = \tilde{\mathbf{c}} \oplus \begin{bmatrix} \tilde{\mathbf{G}} & \tilde{\mathbf{G}}_{\mathbf{d}} \end{bmatrix} B_\infty(\tilde{\mathbf{A}}, \tilde{\mathbf{b}}) \oplus \triangleleft(\mathbf{L}, (\mathbf{c} - \mathbf{h}) \oplus 2\mathbf{G}B_\infty(\mathbf{A}, \mathbf{b}))$ and $\bar{X}_k$ has $(1/2)m_g^2 + (5/2)m_g + 2n$ generators and $(1/2)m_c^2 + (5/2)m_c$ constraints. Therefore,

$$O(\text{CZFO}_{\text{prediction}}) = O(\eta(\mathbf{h}) \oplus \nabla^T \eta(\mathbf{h})(Z - \mathbf{h}) \oplus R) + O(R)$$
$$= O(n) + O(nmm_g) + O(R) = O(nmm_g) + O(R).$$

But

$$O(R) = O(\tilde{\mathbf{c}} \oplus \begin{bmatrix} \tilde{\mathbf{G}} & \tilde{\mathbf{G}}_{\mathbf{d}} \end{bmatrix} B_\infty(\tilde{\mathbf{A}}, \tilde{\mathbf{b}}) \oplus \triangleleft(\mathbf{L}, (\mathbf{c} - \mathbf{h}) \oplus 2\mathbf{G}B_\infty(\mathbf{A}, \mathbf{b})))$$
$$= O(n) + O(\tilde{\mathbf{c}}) + O(\tilde{\mathbf{G}}) + O(\tilde{\mathbf{G}}_{\mathbf{d}}) + O(\tilde{\mathbf{A}}) + O(\tilde{\mathbf{b}}) + O(\mathbf{L}) + O(\triangleleft(\mathbf{L}, (\mathbf{c} - \mathbf{h}) \oplus 2\mathbf{G}B_\infty(\mathbf{A}, \mathbf{b})))$$
$$+ nO(\tilde{\mathbf{Q}}^{[q]}) + O(\text{interval hull of } Z)$$
$$= O(n) + O(nm_g) + O(nm_g^2) + O(nm_g^2) + O(m_g m_c^2 + m_g^2 m_c^2) + O(m_g m_c^2) + O(nm^2)$$
$$+ O(n + nmm_g + m_c(m_g + m_c)^3 + mm_g^2 m_c) + nO(m^2 m_g + mm_g^2) + O(mm_g(m_g + m_c)^3)$$
$$= O(n(m^2 m_g + mm_g^2) + (mm_g + m_c)(m_g + m_c)^3).$$

Therefore,

$$O(\text{CZFO}_{\text{prediction}}) = O(nmm_g) + O(R) = O(n(m^2 m_g + mm_g^2) + (mm_g + m_c)(m_g + m_c)^3).$$

### 5.2. Update step

Let $\bar{X}_k = \{\mathbf{G}_x, \mathbf{c}_x, \mathbf{A}_x, \mathbf{b}_x\} \subset \mathbb{R}^n$, $V = \{\mathbf{G}_v, \mathbf{c}_v, \mathbf{A}_v, \mathbf{b}_v\} \subset \mathbb{R}^{n_v}$, $\mathbf{u}_k \in \mathbb{R}^{n_u}$, and $\mathbf{y}_k \in \mathbb{R}^{n_y}$ with $\mathbf{c}_x \in \mathbb{R}^n$, $\mathbf{G}_x \in \mathbb{R}^{n \times \bar{n}_g}$, $\mathbf{A}_x \in \mathbb{R}^{\bar{n}_c \times \bar{n}_g}$, $\mathbf{b}_x \in \mathbb{R}^{\bar{n}_c}$, $\mathbf{c}_v \in \mathbb{R}^{n_v}$, $\mathbf{G}_v \in \mathbb{R}^{n_v \times n_{g_v}}$, $\mathbf{A}_v \in \mathbb{R}^{n_{c_v} \times n_{g_v}}$, and $\mathbf{b}_w \in \mathbb{R}^{n_{c_v}}$. The update step is performed according to

$$\hat{X}_k = \bar{X}_k \cap_\mathbf{C} ((\mathbf{y}_k - \mathbf{D}_u \mathbf{u}_k) \oplus (-\mathbf{D}_v V)),$$

where $\mathbf{C} \in \mathbb{R}^{n_y \times n}$, $\mathbf{D}_u \in \mathbb{R}^{n_y \times n_u}$, and $\mathbf{D}_v \in \mathbb{R}^{n_y \times n_v}$. As a result, from the prediction step of CZFO, $\bar{X}_k$ has $\bar{n}_g = (1/2)m_g^2 + (5/2)m_g + 2n$ generators and $\bar{n}_c = (1/2)m_c^2 + (5/2)m_c$ constraints, where $m_g = n_g + n_{g_w}$, and $m_c = n_c + n_{c_w}$. Therefore,

$$O(\text{CZFO}_{\text{update}}) = O(\bar{X}_k \cap_\mathbf{C} ((\mathbf{y}_k - \mathbf{D}_u \mathbf{u}_k) \oplus (-\mathbf{D}_v V)))$$
$$= O(n((1/2)m_g^2 + (5/2)m_g + 2n)n_y + n_y n_{g_v}) + O(n_y) + O(n_y + n_y n_u) + O(n_y n_v n_{g_v})$$
$$= O(n_y n(m_g^2 + n) + n_y n_u + n_y n_v n_{g_v}).$$

*5.3. Order reduction*

Let $n_g$ and $n_c$ denote the desired number of generators and constraints of $\hat{X}_k$ after order reduction, i.e., $\hat{X}_k$ must have the same complexity as $\hat{X}_{k-1}$. Order reduction is performed by eliminating $k_c$ constraints and then $k_g$ generators from $\hat{X}_k$. For simplicity, define $m_g = n_g + n_{g_w}$, $m_c = n_c + n_{c_w}$, $\delta_n = n - n_y$, $\delta_v = n_{g_v} - n_{c_v}$, and $\tilde{\delta} = m_g^2 - m_c^2$. The prediction and update steps of CZFO lead to $\hat{X}_k$ with $(1/2)m_g^2 + (5/2)m_g + 2n + n_{g_v}$ generators and $(1/2)m_c^2 + (5/2)m_c + n_{c_v} + n_y$ constraints. Therefore, $k_c = (1/2)m_c^2 + (5/2)m_c + n_{c_v} + n_y - n_c$ and $k_g = (1/2)m_g^2 + (5/2)m_g + 2n + n_{g_v} - n_g - k_c = (1/2)m_g^2 + (5/2)m_g + 2n + n_{g_v} - n_g - (1/2)m_c^2 - (5/2)m_c - n_{c_v} - n_y + n_c$. Consequently,

$$\begin{aligned}
O(\text{eliminate } k_c \text{ constraints from } \hat{X}_k) &= O(k_c(m_g^2 + m_g + n + n_{g_v} + m_c^2 + m_c + n_{c_v} + n_y)^3 \\
&\quad + k_c n(m_g^2 + m_g + n)^2) \\
&= O((m_c^2 + n_{c_v} + n_y)(m_g^2 + n + n_{g_v} + m_c^2 + n_{c_v} + n_y)^3 \\
&\quad + (m_c^2 + n_{c_v} + n_y)n(m_g^2 + n)^2) \\
&= O((m_c^2 + n_{c_v} + n_y)(m_g^2 + m_c^2 + n_{g_v} + n_{c_v} + n + n_y)^3).
\end{aligned}$$

Once again, by eliminating $k_c$ constraints from $\hat{X}_k$ the same number of generators (i.e., $k_c$) is also removed [1]. Hence,

$$\begin{aligned}
O(\text{elim. } k_g \text{ gen. from } \hat{X}_k \text{ after const. elim.}) &= O((n + n_c)^2(m_g^2 + n + n_{g_v} - m_c^2 - n_{c_v} - n_y) \\
&\quad + k_g(m_g^2 + n + n_{g_v} - m_c^2 - n_{c_v} - n_y)(n + n_c)) \\
&= O((n + n_c)^2(\tilde{\delta} + \delta_n + \delta_v) \\
&\quad + (\tilde{\delta} + \delta_n + \delta_v)(\tilde{\delta} + \delta_n + \delta_v)(n + n_c)) \\
&= O((n + n_c)^2(\tilde{\delta} + \delta_n + \delta_v) + (\tilde{\delta} + \delta_n + \delta_v)^2(n + n_c)).
\end{aligned}$$

Finally,

$$\begin{aligned}
O(\text{CZFO}_{\text{reduction}}) &= O(\text{eliminate } k_c \text{ constraints from } \hat{X}_k) + O(\text{elim. } k_g \text{ gen. from } \hat{X}_k \text{ after const. elim.}) \\
&= O((m_c^2 + n_{c_v} + n_y)(m_g^2 + m_c^2 + n_{g_v} + n_{c_v} + n + n_y)^3 \\
&\quad + (n + n_c)^2(\tilde{\delta} + \delta_n + \delta_v) + (\tilde{\delta} + \delta_n + \delta_v)^2(n + n_c)).
\end{aligned}$$




\end{document}